\documentclass[11pt]{article}

\usepackage{latexsym,amssymb,amsmath,natbib,graphicx, amsthm}%,bbm}
\usepackage[margin=1in,left=1in]{geometry}
\usepackage{natbib}
\renewcommand{\cite}{\citeyearpar}

\newcommand {\ME}{\mathbb{E}^{x}}
\newcommand {\MEX}{\mathbb{E}^{X_{T_1}}}
\newcommand {\D}{e^{-\alpha\tau}}

\newcommand {\myBox}{\hspace{\stretch{1}}$\diamondsuit$}
\renewcommand{\S}{\mathcal{S}}

\numberwithin{equation}{section}
\newtheorem{proposition}{Proposition}[section]

\newtheorem{remark}{Remark}[section]
\newtheorem{lemma}{Lemma}[section]
\newtheorem{example}{Example}[section]

\newcommand {\R}{\mathbb{R}}
\newcommand {\F}{\mathcal{F}}

\newcommand {\p}{\mathbb{P}}

\newcommand {\E}{\mathbb{E}}

\title{A Direct Solution Method for Stochastic Impulse Control Problems of One-dimensional Diffusions \footnote{The earlier version of this work was circulated under the title of \emph{Solving Stochastic Impulse Control Problems via Optimal Stopping for One-dimensional Diffusions.}}}
\author{Masahiko Egami \footnote{Department of Mathematics, University
of Michigan, Ann Arbor, MI 48189-1043; egami@umich.edu}}
\date{September 2006
\footnote{ I thank Savas Dayanik for the guidance with respect to this work.  I also thank Erhan Bayraktar for valuable comments and  the participants
at the INFORMS 2004 Annual Meeting in Denver, CO, at the Civitas
Foundation Finance Seminar in Princeton, NJ and at the seminars in various universities.}}

\begin{document}
\maketitle
\begin{abstract}
\noindent We consider stochastic impulse control problems where the
process is driven by a general one-dimensional diffusions. Impulse
control problems are widely used to financial
engineering/decision-making problems such as dividend payout
problem, portfolio optimization with transaction costs, and
inventory control.  We shall show a new mathematical
characterization of the value function as a linear function in
certain transformed space.  Our approach can (1) relieve us from the
burden of guessing and proving the optimal strategy, (2) present a
simple method to find the value function
    and the corresponding control policies, and
(3) handle systematically a broader class of reward and cost
functions than the conventional methods of quasi variational
inequalities, especially because the existence of the finite value
function can be shown in much a simpler way.\\
\\
Key Words: Stochastic impulse control, Diffusions, Optimal stopping, Concavity.\\
AMS Subject Classification: Primary: 49N25  Secondary: 60G40.
\end{abstract}

%---------------------------
% Introduction
%===========================
\section{Introduction}
This paper proposes a general solution method of stochastic impulse
control problems for one dimensional diffusion processes. Stochastic
impulse control problems have attracted a growing interest of many
researchers for the last two decades.  Under a typical setting, the
controller faces some underlying process and reward/cost structure.
There exist continuous and instantaneous components of reward/cost
functions.  By exercising impulse controls, the controller moves the
underlying process from one point to another.  At the same time, the
controller receives  rewards associated with the instantaneous
shifts of the process.  Then the controller's objective is to
maximize the total discounted expected net income.

The mathematical framework to these types of problems is in
Bensoussan and Lions \cite{BL1984}.  Impulse control has been
studied widely in inventory control (Harrison et al.
\cite{HST1983}), exchange rate problem (Jeanblanc-Picqu\'{e}
\cite{MJ1993}, Mundaca and {\O}ksendal \cite{MO1998}, Cadenillas and
Zapatero \cite{CZ2000}), dividend payout problems
(Jeanblac-Picqu\'{e} and Shiryaev \cite{JBS1995}), and portfolio
optimization with transaction costs (Korn \cite{KO1998}, Morton and
Pliska \cite{MP1995}).  Korn \cite{KO1999} surveys the applications
in mathematical finance. Also see Chancelier et al.\cite{COS2002}
for a combination of optimal stopping and impulse control problems.
In many economic and financial applications where the controlled
process is described as an It\^{o} diffusion, the solution to the
problem demands a through study of a related Hamilton-Jacobi-Bellman
equation and quasi-variational inequalities.  The method of
quasi-variational inequalities split by a guess the state space into
intervention and no intervention (continuation) regions. One guesses
the form of (a) continuation region, (b) associated optimal policy,
and (c) the value function.  Then optimality of the candidate policy
must be verified.  Both steps are often very difficult and the
success depends heavily on the form of the controlled process,
reward and cost functions.

Alternatively, an impulse control problem can be viewed as a
sequence of optimal stopping problems.  The connection between
impulse control and optimal stopping has been investigated by Davis
\cite{DV1992} and {\O}ksendal and Sulem \cite{OS2002} among others.
In this setting, the value functions of a sequence of optimal
stopping problems converge to the value function of the impulse
control problem under suitable conditions.

In this paper, we utilize this connection together with a novel
method of Dayanik and Karatzas~\cite{DK2003} for optimal stopping
problems.  We use it to identify a new and useful characterization
of the solution of the original impulse control problem.  At the end
we get rid of the sequence of optimal stopping problems altogether:
the new characterization allow us to propose a new direct solution
method for impulse control problems.

In the next section, we briefly go over the solution method for
optimal stopping problems of one-dimensional diffusions.  We
describe the impulse control problem and its solution in Section 3.
Examples are presented in Section 4.  Finally, extensions and
concluding remarks are in Section 5.

%-----------------------
% Section 2
%=======================

\section{Summary of the Key Results of Optimal Stopping}
\label{sec:summary} Let $(\Omega, \F, \p)$ be a complete
probability space with a standard Brownian motion $W=\{W_t; t\geq
0\}$ and consider the diffusion process $X^0$ with state pace
$\mathcal{I}\subseteq \mathbb{R}$ and dynamics
\begin{equation}\label{eq:process}
dX^0_t=\mu(X^0_t)dt + \sigma(X^0_t)dW_t
\end{equation}
for some Borel functions $\mu :\mathcal{I}\rightarrow \mathbb{R}$
and $\sigma :\mathcal{I}\rightarrow (0, \infty)$.  We emphasize
here that $X^0$ is an uncontrolled process.  We assume that
$\mathcal{I}$ is an interval with endpoints $-\infty\leq a < b
\leq+\infty$, and that $X^0$ is regular in $(a, b)$; in other
words, $X^0$ reaches $y$ with positive probability starting at $x$
for every $x$ and $y$ in $(a,b)$.  We shall denote by
$\mathbb{F}=\{\mathcal{F}_t\}$ the natural filtration generated by
$X^0$.

Let $\alpha \geq 0$ be a real constant and $h(\cdot)$ a Borel
function such that $\ME[e^{-\alpha \tau}h(X^0_{\tau})]$ is
well-defined for every $\mathbb{F}$-stopping time $\tau$ and $x\in
\mathcal{I}$.  Let $\tau_y$ be the first hitting time of $y\in
\mathcal{I}$ by $X^0$, and let $c\in\mathcal{I}$ be a fixed point
of the state space.  We set:
\begin{align} \nonumber
\begin{aligned}
    \psi(x) &=  \begin{cases}
                 \ME[e^{-\alpha\tau_c}1_{\{\tau_c<\infty\}}], & x\leq c, \\
                 1/\E^{c}[e^{-\alpha\tau_x}1_{\{\tau_x<\infty\}}],
                 &x>c,\end{cases}
\hspace{0.4cm}
    \varphi(x) &= \begin{cases}
                 1/\E^{c}\left[e^{-\alpha\tau_x}1_{\{\tau_x<\infty\}}\right], & x\leq c, \\
                 \ME[e^{-\alpha\tau_c}1_{\{\tau_c<\infty\}}],
                 &x>c,
     \end{cases}
\end{aligned}
\end{align}
and
\begin{align} \label{eq:F}
F(x)&\triangleq\frac{\psi(x)}{\varphi(x)}, \hspace{0.5cm} x\in
\mathcal{I}.
\end{align}
Then $F(\cdot)$ is continuous and strictly increasing.  It should
be noted that $\psi(\cdot)$ and $\varphi(\cdot)$ consist of an
increasing and a decreasing solution of the second-order
differential equation $(\mathcal{A}-\alpha)u=0$ in $\mathcal{I}$
where $\mathcal{A}$ is the infinitesimal generator of $X^0$.  They
are linearly independent positive solutions and uniquely
determined up to multiplication.  For the complete
characterization of $\psi(\cdot)$ and $\varphi(\cdot)$
corresponding to various types of boundary behavior, refer to
It\^{o} and McKean \cite{IM1974}.

Let $F :[c, d]\rightarrow\mathbb{R}$ be a strictly increasing
function.  A real valued function $u$ is called \emph{$F$-concave}
on $[c, d]$ if, for every $a\leq l<r\leq b$ and $x\in[l, r]$,
\begin{equation*}
u(x)\geq
u(l)\frac{F(r)-F(x)}{F(r)-F(l)}+u(r)\frac{F(x)-F(l)}{F(r)-F(l)}.
\end{equation*}
We denote by
\begin{equation}\label{eq:value}
V(x)\triangleq \sup_{\tau\in\S}\ME[\D h(X^0_\tau)], \hspace{0.5cm}
x\in[c, d]
\end{equation}
the value function of the optimal stopping problem with the reward
function $h(\cdot)$ where the supremum is taken over the class
$\S$ of all $\mathbb{F}$-stopping times.  Then we have the
following results, the proofs of which we refer to Dayanik and
Karatzas~\cite{DK2003}.
%\begin{proposition}\normalfont\label{prop:1}
%For a given function $U$: $[c, d]\rightarrow[0,+\infty)$ the
%quotient $U(\cdot)/\varphi(\cdot)$ is an $F$-concave function if
%and only if $U(\cdot)$ is $\alpha$-excessive, i.e.,
%\begin{align}
%U(x)\geq \ME[e^{-\alpha\tau}U(X^0_\tau)], &\forall \tau\in\S,
%\forall x\in[c,d].
%\end{align}
%\end{proposition}
\begin{proposition}\normalfont \label{prop:2}
The value function $V(\cdot)$ of (\ref{eq:value}) is the smallest
nonnegative majorant of $h(\cdot)$ such that
$V(\cdot)/\varphi(\cdot)$ is $F$-concave on $[c,d]$.
\end{proposition}
\begin{proposition}\normalfont \label{prop:3}
Let $W(\cdot)$ be the smallest nonnegative concave majorant of
$H\triangleq (h/\varphi)\circ F^{-1}$ on $[F(c), F(d)]$, where
$F^{-1}(\cdot)$ is the inverse of the strictly increasing function
$F(\cdot)$ in (\ref{eq:F}).  Then $V(x)=\varphi(x)W(F(x))$ for
every $x\in[c, d]$.
\end{proposition}
\begin{proposition}\normalfont \label{prop:4}
Define
\begin{equation}\label{eq:opt}
S\triangleq\{x\in[c, d]: V(x)=h(x)\}, \hspace{0.5cm}\text{and}
\hspace{0.5cm} \tau^{*}\triangleq \inf\{t\geqq0: X^0_t\in S\}.
\end{equation}
If $h(\cdot)$ is continuous on $[c, d]$, then $\tau^{*}$ is an
optimal stopping rule.
\end{proposition}

\section{Impulse control problems and its solution} \label{sec:abp}
Suppose that at any time $t\in\mathbb{R_+}$ and any state
$x\in\mathbb{R_+}$, we can intervene and give the system an impulse
$\xi\in \mathbb{R}$. Once the system gets intervened, the point
moves from $x$ to $y\in\mathbb{R_+}$ with associated rewards earned.
An impulse control for the system is a double sequence,
\begin{equation*}
\nu =(T_1, T_2,....T_i....; \xi_1,\xi_2,...\xi_i....)
\end{equation*}
where $0\leq T_1<T_2<....$ are an increasing sequence of
$\mathbb{F}$-stopping times and $\xi_1$, $\xi_2...$ are
$\mathcal{F}_{T_i}$-measurable random variables representing
impulses exercised at the corresponding intervention times $T_i$
with $\xi_i\in Z$ for all $i$ where $Z\in\mathbb{R}$ is a given
set of admissible impulse values.  The controlled process is, in
general, described as follows:
\begin{eqnarray} \label{eq:control}
dX_t &=& \mu(X_t)dt + \sigma(X_t)dW_t, \quad T_{i-1}\leq t<T_i \\
X_{T_{i}}&=&\Gamma(X_{T_{i}-}, \xi_i)
\end{eqnarray}
with some mapping $\Gamma: \mathbb{R}\times\mathbb{R}\rightarrow
\mathbb{R}$.

In this section, we consider the absorbing boundary problem.  Let
$0$ be the absorbing state, without loss of generality, and $\tau_0
\triangleq \inf\{t: X_t=0 \}$ the ruin time.  With the absorbing
state at $0$, it is natural to consider a set of problems where
$Z\in \mathbb{R_+}$ (i.e., $\xi_i=x_i-y_i>0$ for all $i$) and
$X_{T_i}=X_{T_i-}-\xi_i$.  (We shall comment on cases where
interventions are allowed in both positive and negative directions
in section \ref{sec:conclusion}.)

With each pair $(T_i, \xi_i)$, we associate
the interventions  %$Z$ on \{$t<\tau_0$\} defined by
\begin{equation}\label{eq:Z}
K(X_{T_i-},X_{T_i})
\end{equation}
where $K(x, y):\mathbb{R}\times\mathbb{R}\rightarrow \mathbb{R}$
is a given continuous function in the first and second argument
that represents benefit/cost at interventions. Our result below
does not depend on the specification of $K(\cdot)$.
% In many cases, for example, $K(\cdot)$ is in the form,
%\begin{align*}\label{eq:K}
%K(X_{T_i-}, X_{T_i})&=h(X_{T_i-},X_{T_i})-c(X_{T_i-}, X_{T_i})
%\end{align*}
%where we call $h$ and $c$, some Borel functions $h,
%c:\mathbb{R}\times\mathbb{R}\rightarrow \mathbb{R}$ reward and
%cost function, respectively.
We assume that, for any point $x\in\mathbb{R}$,
\begin{equation}\label{eq:Kneg}
K(x, x)<0.
\end{equation}
due to the fixed cost incurred.  We consider the following
performance measure with $\nu\in\mathcal{V}$, a collection of
admissible strategies,
\begin{equation}\label{eq:J}
    J^\nu(x)=\ME\left[\int_0^{\tau_0}e^{-\alpha s}f(X_s)ds+e^{-\alpha\tau_0}P+\sum_{T_i<\tau_0}e^{-\alpha T_i}K(X_{T_i-},X_{T_i})\right]
\end{equation}
where $P\in \R_{-}$ is a constant penalty at the ruin time and
$f:\mathbb{R}\rightarrow \mathbb{R}$ is a continuous function,
satisfying :
\begin{equation} \label{eq:fcon}
\ME\left[\int_0^\infty e^{-\alpha s}|f(X_s)|ds\right]<\infty.
\end{equation}
Our goal is to find the optimal strategy $\nu^*(T_i,
\xi_i)_{i\geq0}$ and the corresponding value function,
\begin{equation} \label{eq:impulsevalue}
v(x)\triangleq\sup_{\nu\in\mathcal{V}} J^{\nu}(x)=J^{\nu^*}(x).
\end{equation}
Let us briefly go over our plan.  In section \ref{sec:first} we
shall characterize optimal intervention times $T_i$ as exit times of
the process $X$ from an interval by implementing recursive optimal
stopping scheme that eventually solves the original impulse control
problem. Using the results, in section \ref{sec:sec}, we consider a
special case where the mapping $x\rightarrow\frac{K}{\varphi}(x)$
$:\mathbb{R_+}\rightarrow\mathbb{R_+}$ is $F$-concave.  We show,
under this assumption, that the optimal intervention times $T_i$ are
characterized as exit times from an interval, say $(0, b^*)$ for
every $i$.  Then we characterize the value function for impulse
control problems and present a solution method based on the
characterization of the intervention times and value function. In
section \ref{sec:gen}, we consider the general case where the
$F$-concavity assumption above does not hold.

%-------------------------------------
% 3.1 Characterization
%====================================
\subsection{A sequence of optimal stopping problems} \label{sec:first}
In this subsection, we consider a recursive optimal stopping with a
view to characterizing intervention times for the impulse control
problems.  Here we assume that no absorbing boundary exists.  As we
will see in the next subsection, the existence of an absorbing state
is easily incorporated.  Hence by using the same $v(x)$, we consider
the problem,
\begin{equation}\label{eq:newv}
v(x)=\sup_{\nu}\ME\left[\int_0^{\infty}e^{-\alpha
s}f(X_s)ds+\sum_{i}e^{-\alpha T_i}K(X_{T_i-},X_{T_i})\right],
\end{equation}
and define the set $S_n$ and the objective function $v_n$ as
follows:
\begin{equation*}
    S_n\triangleq \{\nu\in S; \nu=(T_1, T_2,... T_{n+1}; \xi_1,
    \xi_2,...\xi_n); T_{n+1}=+\infty\},
\end{equation*}
and
\begin{equation}\label{eq:vn}
v_n(x)\triangleq \sup_{\nu\in
S_n}\ME\left[\int_0^{\infty}e^{-\alpha
s}f(X_s)ds+\sum_{T_i}e^{-\alpha T_i}K(X_{T_i-},X_{T_i})\right].
\end{equation}
%Define $\hat{g}(x)\triangleq \ME[\int_0^{\infty}e^{-\alpha
%s}f(X_s)ds]$. Note that this $X$ is the controlled process and
%$g(x)\neq\hat{g}(x)$
In other words, we are allowed to make at most $n$ interventions.
For this recursive approach, see, for example, Davis \cite{DV1992}
and {\O}ksendal and Sulem \cite{OS2002}.  We use the following
simple notation:
\begin{equation} \label{eq:gf}
g(x)\triangleq\ME\left[\int_0^\infty e^{-\alpha
s}f(X_s^0)ds\right].
\end{equation}
Let $\mathcal{H}$ denote the space of all Borel functions. Define
the two operators $\mathcal{M} :
\mathcal{H}\rightarrow\mathcal{H}$ and $\mathcal{L}
:\mathcal{H}\rightarrow\mathcal{H}$ as follows:
\begin{equation}\label{eq:M}
\mathcal{M}u(x)\triangleq\sup_{y\in \mathbb{R}}[K(x,
y)-(g(x)-g(y))+u(y)],
\end{equation}
and
\begin{equation}\label{eq:L}
\mathcal{L}u(x)\triangleq\sup_{\tau\in\mathcal{S}}\ME[e^{-\alpha
\tau}\mathcal{M}u(X_{\tau-})],
\end{equation}
for $u\in\mathcal{H}$.  From the definition of the two operators,
$a_1(x)\leq a_2(x)$ for $x\in\mathbb{R}, a_1(\cdot),
a_2(\cdot)\in\mathcal{H}$ implies $\mathcal{M}a_1(x)\leq
\mathcal{M}a_2(x)$ and $\mathcal{L}a_1(x)\leq \mathcal{L}a_2(x)$ for
all $x\in\mathbb{R}$. Consider the following recursive formula:
\begin{equation} \label{eq:wrec}
w_{n+1}(x)=\sup_{\tau\in\mathcal{S}, \xi}\ME\left[\int_0^\tau
e^{-\alpha s}f(X_s)ds + e^{-\alpha\tau}(K(X_{\tau-},
X_{\tau})+w_n(X_\tau))\right],
\end{equation}
which is equivalent to
\begin{equation} \label{eq:wrec2}
w_{n+1}(x)-g(x)=\sup_{\tau\in\mathcal{S},
\xi}\ME[e^{-\alpha\tau}(K(X_{\tau-},
X_{\tau})-g(X_{\tau-})+w_n(X_\tau))]
\end{equation}
by applying the strong Markov property with (\ref{eq:fcon}) to the
integral term.  In fact, this derivation is explained in detail in
subsection \ref{sec:sec}.  By defining
\begin{equation*}
\phi\triangleq w-g,
\end{equation*}
and adding and subtracting $g(X_{\tau})$ on the right hand side of
(\ref{eq:wrec2}), it becomes
\begin{equation*}
\phi_{n+1}(x)=\sup_{\tau\in\mathcal{S}}\ME[e^{-\alpha\tau}\mathcal{M}\phi_{n}(X_{\tau-})].
\end{equation*}
It should be noted that, for each $n$, this is an optimal stopping
problem over $\tau$ and can be written, by using the operator
defined in (\ref{eq:L}),
\begin{equation} \label{eq:stops}
\phi_{n+1}(x)=\mathcal{L}\phi_{n}(x).
\end{equation}

Let us start this recursive scheme with $w_0(x)\triangleq g(x)$
(i.e., no interventions are allowed, equivalently $\phi_0(x)=0$) and
define recursively $\phi_n(x)\triangleq
w_n(x)-g(x)=\mathcal{L}(w_{n-1}(x)-g(x))=\mathcal{L}\phi_{n-1}$.
Clearly,
\begin{align*}\label{eq:phi1}
\phi_1(x)&=\mathcal{L}\phi_0(x) =\sup_{\tau\in \S}\ME[e^{-\alpha\tau}(\mathcal{M}(w_0(X_\tau)-g(X_\tau))]\nonumber\\
&=\sup_{\tau \in \S, \xi\in
\R_+}\ME[e^{-\alpha\tau}\left\{K(X_{\tau-},
X_\tau)-g(X_{\tau-})+g(X_\tau)\right\}].
\end{align*}
On the other hand,
\begin{align*}
&v_1(x)-g(x)=\sup_{\nu\in S_1}\ME\left[\int_0^{\infty}e^{-\alpha
s}f(X_s)ds+e^{-\alpha
\tau}K(X_{\tau-},X_{\tau})\right]-\ME\left[\int_0^{\infty}e^{-\alpha
s}f(X_s^0)ds\right]\\
&=\sup_{\tau \in \S, \xi\in \R_+}\ME\Big[\int_0^{\tau}e^{-\alpha
s}f(X_s)ds+e^{-\alpha \tau}\left(\E^{X_{\tau}}\left[\int_0^\infty
e^{-\alpha s}f(X_s)ds\right]+K(X_{\tau-},X_{\tau})\right) \\
&\hspace{2.7cm}-\int_0^\tau e^{-\alpha s}
f(X_s^0)ds-e^{-\alpha\tau}g(X_{\tau-})\Big]\\
&=\sup_{\tau \in \S, \xi\in
\R_+}\ME[e^{-\alpha\tau}\left\{K(X_{\tau-},X_{\tau})+g(X_{\tau})-g(X_{\tau-})\right\}].
\end{align*}
The last equation is due to the fact that only one intervention is
allowed.  Hence we have $w_1(x)=v_1(x)$. By the definition of the
recursive scheme, $w_n$ is an increasing sequence (i.e, $w_1(x)\leq
w_2(x) \leq ....$ for all $x\in\mathbb{R}$).  In fact, we shall
prove that $w_n=v_n$ for all $n$ in Lemma \ref{lem:Davis}.  Before
that, we need the following lemma to relate this recursive scheme
with the
method described in Section \ref{sec:summary}.\\

%---------------------------
% Lemma F-concavity of L
%===========================
\begin{lemma}\label{lem:Lop}
The mapping $x\rightarrow\frac{\mathcal{L}\phi(x)}{\varphi(x)}$
$:\mathbb{R}_+\rightarrow\mathbb{R}_+$ is $F$-concave.
\end{lemma}
\begin{proof}
We shall fix some $x\in(l, r)\subseteq [c, d]$. Since
$\mathcal{M}\phi(\cdot)$ is bounded there, for a given
$\varepsilon
>0$, there are admissible $\varepsilon$-optimal intervention pairs
$(\sigma^l_\varepsilon, \xi^l_\varepsilon)$ and
$(\sigma^r_\varepsilon, \xi^r_\varepsilon)$ such that
\begin{equation*}
  \E^l [e^{-\alpha\sigma^l_\varepsilon} \mathcal{M}\phi(X_{\sigma^l_\varepsilon}) ]> \mathcal{L}\phi(l)-\varepsilon,
  \quad\text{and}\quad
 \E^r [e^{-\alpha\sigma^r_\varepsilon} \mathcal{M}\phi(X_{\sigma^r_\varepsilon}) ]>
 \mathcal{L}\phi(r)-\varepsilon.
\end{equation*}
Define another stopping time $\sigma^{lr}_\varepsilon \in
\mathcal{S}$ with
\begin{eqnarray} \nonumber
    \sigma^{lr}_\varepsilon &\triangleq&  \begin{cases}
                 \tau^l + \sigma^l_\varepsilon\circ \theta_{\tau^l}, &\text{if}\quad \tau^l<\tau^r, \\
                 \tau^r + \sigma^r_\varepsilon\circ \theta_{\tau^r}, &\text{if}\quad \tau^l>\tau^r.
     \end{cases}
\end{eqnarray}
Putting all together, with the strong Markov property of $X$, we
have
\begin{align*}
\mathcal{L}\phi(x)&\geq
\ME[e^{-\alpha\sigma^{lr}_\varepsilon}\mathcal{M}\phi(X_{\sigma^{lr}_\varepsilon})]\\
&>(\mathcal{L}\phi(l)-\varepsilon)\ME[e^{-\alpha\tau^l}1_{\{\tau^l<\tau^r\}}]
+(\mathcal{L}\phi(r)-\varepsilon)\ME[e^{-\alpha\tau^r}1_{\{\tau^l>\tau^r\}}]\\
&\geq
\frac{\mathcal{L}\phi(l)}{\varphi(l)}\varphi(x)\frac{F(r)-F(x)}{F(r)-F(l)}
+\frac{\mathcal{L}\phi(r)}{\varphi(r)}\varphi(x)\frac{F(x)-F(l)}{F(r)-F(l)}-\varepsilon.
\end{align*}
Since $\varepsilon$ is arbitrary, we have an $F$-concavity.
\end{proof}
This lemma guarantees that we can use Proposition \ref{prop:2} to
\ref{prop:4} to identify the value function and an optimal stopping
rule for each of the recursive optimal stopping problems
(\ref{eq:wrec}). Let us define, for notational convenience,
\begin{equation}\label{eq:Kbar}
\bar{K}(x, y)\triangleq K(x, y)-(g(x)-g(y)).
\end{equation}
Further, we prove the following properties of the recursive
optimization scheme.

\begin{lemma} \label{lem:Davis}
If we define $w_n$ by (\ref{eq:wrec}) (with $w_0=g$) and $v_n$ by
(\ref{eq:vn}), then
\begin{gather*}
w_n(x)=v_n(x) \quad \text{for each $n$} \quad\text{and}\quad
v(x)=\lim_{n\rightarrow \infty}w_n(x).
\end{gather*}
Moreover, $w$ is the smallest solution majorizing $g$ of the
functional equation $w-g=\mathcal{L}(w-g)$.
\end{lemma}
\begin{proof}
The proof is given in Appendix.
\end{proof}
Hence if we solve the optimal stopping problem
\begin{equation}\label{eq:recursion}
\phi_{n+1}(x)=\sup_{\tau\in\mathcal{S}}\ME[e^{-\alpha\tau}\mathcal{M}\phi_n(X_{\tau-})]
\end{equation}
recursively for each $n$, then we obtain
$\phi(x)=\lim_{n\rightarrow\infty}\phi_n(x)=\lim_{n\rightarrow\infty}v_n(x)-g(x)=v(x)-g(x)$.
Summarizing the above argument, we have the following proposition:

%--------------------------------------
%  Proposition 3-1
%======================================
\begin{proposition}\label{Fconc}
The value function $v(x)$ for (\ref{eq:newv}) is given by the
smallest solution majorizing $g$ of the functional equation
$v-g=\mathcal{L}(v-g)$, and
$\frac{v-g}{\varphi}(=\frac{\phi}{\varphi})$ is always
$F$-concave.
\end{proposition}
\begin{proof}
The first statement comes from Lemma \ref{lem:Davis}.  By the
recursive method that we described above, we are solving a series
of optimal stopping problems for each $\phi_n$.  Hence Lemma
\ref{lem:Lop} and Proposition \ref{prop:2} give the second
statement.
\end{proof}

%---------------------------------
% 3.2 F-Concave Reward Case
%=================================
\subsection{Characterization of the Intervention Times and the Value Function:
$F$-Concave Reward Case}\label{sec:sec} Based on the results in
the previous subsection, we first consider a special case where
the mapping $x\rightarrow\frac{\bar{K}}{\varphi}(x)$
$:\mathbb{R_+}\rightarrow\mathbb{R_+}$ is $F$-concave.  The
argument in the previous subsection is modified to incorporate the
existence of the ruin state. Instead of (\ref{eq:vn}) and
(\ref{eq:wrec}), we define, respectively,
\begin{align*}
v_n(x)&\triangleq \sup_{\nu\in
S_n}\ME\left[\int_0^{\tau_0}e^{-\alpha
s}f(X_s)ds+e^{-\alpha\tau_0}P%1_{\{\tau_0<T_n\}}
+\sum_{T_i<\tau_0}e^{-\alpha T_i}K(X_{T_i-},X_{T_i})\right]\\
w_{n+1}(x)&\triangleq\sup_{\tau\in\mathcal{S},
\xi}\ME\Big[\int_0^{\tau_0\wedge\tau} e^{-\alpha s}f(X_s)ds
+e^{-\alpha\tau_0}P1_{\{\tau_0<\tau\}} \\
&\hspace{6cm}+ e^{-\alpha\tau}\{K(X_{\tau-},
X_{\tau})+w_n(X_\tau)\}1_{\{\tau<\tau_0\}}\Big]
\end{align*}
with
\begin{equation*}
w_0(x)=\ME\left[\int_0^\infty e^{-\alpha
s}f(X^0_s)1_{\{s<\tau_0\}}ds+e^{-\alpha\tau_0}P\right]\triangleq
g_0(x).
\end{equation*}
Then by defining the operator $\mathcal{L}: \mathcal{H}\rightarrow
\mathcal{H}$ instead of (\ref{eq:L}),
\begin{equation}\label{eq:L2}
\mathcal{L}u(x)\triangleq\sup_{\tau\in\mathcal{S}}\ME[e^{-\alpha
\tau}\mathcal{M}u(X_{\tau-})1_{\{\tau<\tau_0\}}+e^{-\alpha\tau_0}(P-g(0))1_{\{\tau_0<\tau\}}],
\end{equation}
we have the same recursion formula as in (\ref{eq:stops}).  We can
obtain the same results as in Lemma \ref{lem:Lop} and Lemma
\ref{lem:Davis}. Proposition \ref{Fconc} also holds with one change
that the value function is given by the smallest solution majorizing
$g_0$ of the functional equation $v-g=\mathcal{L}(v-g)$ where
$\mathcal{L}: \mathcal{H}\rightarrow \mathcal{H}$ is given by
(\ref{eq:L2}). Now we consider the characterization of the
intervention times.
\begin{proposition}\label{prop:exitchar}
If the mapping $x\rightarrow\frac{\bar{K}}{\varphi}(x)$
$:\mathbb{R_+}\rightarrow\mathbb{R_+}$ is $F$-concave and 0 is an
absorbing state, then the optimal intervention times $T_i^*$ are
given, for some $b^*\in\mathbb{R_+}$, by
\begin{equation*}
T_i^*=\inf\{t>T_{i-1}^*; X_t\notin(0, b^*), \quad i=1,2,....\}.
\end{equation*}
\end{proposition}

\begin{proof}
Our proof is constructive, describing the procedure of recursive
optimization steps.  For any $n\geq 1$, in view of Lemma
\ref{lem:Lop}, $\phi_n(x)$ is the smallest $F$-concave majorant of
$\mathcal{L}\phi(x)/\varphi(x)$.  This majorant (that passes $(F(0),
\frac{P-g(0)}{\varphi(0)}$ in the transformed space) always exists.
Indeed, since we consider the case of $\xi_i>0$, i.e, $x> y$ for
$K(x, y)$ and
\begin{align*}
\mathcal{M}\phi_0(x)&=\sup_{y\in \mathbb{R}_+}[K(x,
y)-(g(x)-g(y))+\phi_0(y)]=\sup_{y\in \mathbb{R}_+}[K(x,
y)-(g(x)-g(y))+(g_0(y)-g(y))],
\end{align*}
we should check whether the concave majorant exists, namely,
\begin{equation}\label{eq:check0}
\lim_{x\downarrow 0}(K(x, y)-g(x)+g_0(y))<P-g(0)
\end{equation}
holds when $y\downarrow 0$. Note that $\lim_{y\downarrow 0}g_0(y)=
P$ and $g(x) \rightarrow g(0)$ as $x\rightarrow 0$ due to the
continuity of $f$. Hence (\ref{eq:check0}) holds in the neighborhood
of $y=0$ because of (\ref{eq:Kneg}). In the subsequent iterations,
we consider
\begin{align*}
\mathcal{M}\phi_1(x)&=\sup_{y\in \mathbb{R}_+}[K(x,
y)-(g(x)-g(y))+\phi_1(y)].
%&=\sup_y[K(x, y)-(g(x)-g(y))+(w_1(y)-g(y))]\\
%&=\sup_y[K(x, y)-g(x)+w_1(y)]
\end{align*}
We should check if the expression inside the supremum operator
becomes less than $P-g(0)$ as $x\downarrow 0$ and $y\downarrow 0$.
Since $\lim_{y\downarrow 0}\phi_1(y)=\phi_1(0)=P-g(0)$ by the
concavity (hence continuity) of $\phi_1$ and since
$\lim_{x\downarrow 0}g(x)=\lim_{y\downarrow 0}g(y)$, we have in the
neighborhood of $y=0$,
\begin{align*}
\lim_{x\downarrow 0}K(x, 0)+P-g(0)<P-g(0)
\end{align*}
holds.  Hence the concave majorant always exist also in the
subsequent iterations.

Now the $F$-concavity of $\phi_n$ is obviously maintained for all
$n$. The limit function,
$\phi(x)\triangleq\lim_{n\rightarrow\infty}\phi_n(x)$ exists and is
also $F$-concave.  Accordingly, $\bar{K}(x, y)/\varphi(x)+\phi(y)$
is $F$- concave. Hence $\phi$ and $\bar{K}+\phi$ meet once and only
once. Recall that the value function satisfies
$\phi=\mathcal{L}\phi$. This implies that the continuous region is
in the form of $(0, b^*)$ for some
$b^*\in\mathbb{R_+}$, which completes the proof. %Moreover, if
%$K(x, y)$ is $F$-concave in $x$, then there are no two (or more) $y^i$ and $y^j$ such that
%\begin{equation*}
%K(x, y^i)=K(x, y^j) \quad \forall x\in\mathbb{R_+}.
%\end{equation*}
%Hence $b^*$ is unique.
\end{proof}

%\begin{remark} In the above proof, note that for the limit function $\phi$, we have
%\begin{eqnarray*} \nonumber
%   K(x, y^*)+\phi(y^*)&<& \phi(x), \quad 0 \leq x < b^*, \\
%   K(x, y^*)+\phi(y^*)&=& \phi(x), \quad x\geq b^*.
%\end{eqnarray*}
%by the construction. \myBox
%\end{remark}

By using the above characterization of intervention times, we next
want to characterize the value function and reduce the impulse
control problem (\ref{eq:impulsevalue}) to some optimal stopping
problem.  Moreover, we shall present a method that does not have
to go through the iteration scheme.  Let us first simplify
$J^\nu$:
\begin{equation} \label{eq:u}
J^\nu(x)=\ME\left[\int_0^{\tau_0}e^{-\alpha
s}f(X_s)ds+e^{-\alpha\tau_0}P+\sum_{T_i<\tau_0}e^{-\alpha
T_i}K(X_{T_i-},X_{T_i})\right].
\end{equation}
This is just a reproduction of (\ref{eq:J}).  Let us split the
right hand side of (\ref{eq:u}) into pieces and use the strong
Markov property (together with the shift operator $\theta(\cdot)$)
to each of them.  The first term becomes
\begin{align*}
&\ME\left[\int_0^{\tau_0}e^{-\alpha
s}f(X_s)ds\right]=\ME\left[1_{\{T_1<\tau_0\}}\left\{\int_0^{T_1}e^{-\alpha
s}f(X_s)ds+ e^{-\alpha T_1}\MEX \int_{0}^{\infty}e^{-\alpha
s}f(X_s)1_{\{s<\tau_0\}}ds\right\}\right]\\
&\hspace{7cm}+\ME\left[1_{\{T_1>\tau_0\}}\int_0^{\tau_0}f(X_s)ds\right]
\end{align*}The second and third terms become
\begin{align*}
\ME e^{-\alpha\tau_0}P &=\ME\left[1_{\{T_1<\tau_0\}}e^{-\alpha
T_1}\ME[e^{-\alpha(\tau_0-T_1)}P|\mathcal{F}_{T_1}]\right]+\ME\left[1_{\{T_1>\tau_0\}}e^{-\alpha\tau_0}P\right]\\
&=\ME\left[1_{\{T_1<\tau_0\}}e^{-\alpha
T_1}\ME[e^{-\alpha(\tau_0\circ\theta(T_1))}P|\mathcal{F}_{T_1}]\right]+\ME\left[1_{\{T_1>\tau_0\}}e^{-\alpha\tau_0}P\right]\\
&=\ME\left[1_{\{T_1<\tau_0\}}e^{-\alpha T_1}\MEX
(e^{-\alpha\tau_0}P)\right]+\ME\left[1_{\{T_1>\tau_0\}}e^{-\alpha\tau_0}P\right]
\end{align*}
and
\begin{align*} &\ME\left[\sum_{T_i<\tau_0}e^{-\alpha
T_i}K(X_{T_i-},X_{T_i})\right]\\
&= \ME\left[1_{\{T_1<\tau_0\}}\left\{e^{-\alpha
T_1}K(X_{T_1-},X_{T_1})+e^{-\alpha T_1}\sum_{i=2}e^{-\alpha
(T_i-T_1)}K(X_{T_i-},X_{T_i})1_{\{T_i<\tau_0\}}\right\}\right]\nonumber \\
%&=\ME\left[1_{\{T_1<\tau_0\}}\left\{e^{-\alpha
%T_1}K(X_{T_1-},X_{T_1}) +e^{-\alpha
%T_1}\ME\left[\sum_{i=2}e^{-\alpha
%(T_{i-1}\circ\theta(T_1))}K(X_{S_{(i-1)-}},X_{S_{(i-1)}})1_{\{T_i<\tau_0\}}|\mathcal{F}_{T_1}\right]\right\}\right] \nonumber\\
&=\ME\left[1_{\{T_1<\tau_0\}}\left\{e^{-\alpha
T_1}K(X_{T_1-},X_{T_1})+e^{-\alpha
T_1}\ME[\sum_{T_i<\tau_0}e^{-\alpha
(T_i\circ\theta(T_1))}K(X_{S_{i-}},X_{S_i})|\mathcal{F}_{T_1}]\right\}\right] \nonumber\\
&=\ME\left[1_{\{T_1<\tau_0\}}e^{-\alpha
T_1}\left\{K(X_{T_1-},X_{T_1})+\MEX\sum_{i=1}e^{-\alpha
T_i}K(X_{T_i-},X_{T_i})1_{\{T_i<\tau_0\}}\right\}\right],
\end{align*}
where $S_i\triangleq T_1+T_{i}\circ\theta(T_1)$ and the index $i$
runs from $1$ for the sum in the second equality. Combining the
three terms and rearranging, we have
\begin{multline} \label{eq:simple}
J^\nu(x)=\ME\left[1_{\{T_1<\tau_0\}}\left\{\int_0^{T_1}e^{-\alpha
s}f(X_s)ds+e^{-\alpha T_1}K(X_{T_1-},X_{T_1})+e^{-\alpha
T_1}J^\nu(X_{T_1})\right\}\right]\\
+\ME\left[1_{\{T_1>\tau_0\}}\left\{\int_0^{\tau_0}e^{-\alpha
s}f(X_s)ds+e^{-\alpha\tau_0}P\right\}\right].
\end{multline}
For any $\mathbb{F}$ stopping time $\tau$, the strong Markov
property with our assumption (\ref{eq:fcon}) gives us
\begin{equation*}
  \ME\left[\int_0^\tau e^{-\alpha s}f(X_s^0)ds\right]=g(x)-\ME\left[ e^{-\alpha
  \tau}g(X_{\tau}^0)\right]
\end{equation*}
where $g(\cdot)$ is defined as in (\ref{eq:gf}).  We apply this
result to (\ref{eq:simple}) by reading $\tau=T_1$ and $\tau=\tau_0$
to derive%. Since
%\begin{equation*}
%\ME\left[\int_0^{T_1} e^{-\alpha
%s}f(X_s)ds\right]=\ME\left[\int_0^{T_1} e^{-\alpha
%s}f(X_s^0)ds\right]=g(x)-\ME\left[e^{-\alpha T_1}g(X^0_{T_1})\right]
%\end{equation*}
%and also, on $\{\tau_0<T_1\}$,
%\begin{equation*}
%\ME\left[\int_0^{\tau_0} e^{-\alpha
%s}f(X_s)ds\right]=\ME\left[\int_0^{\tau_0} e^{-\alpha
%s}f(X_s^0)ds\right]=g(x)-\ME\left[ e^{-\alpha
%\tau_0}g(X^0_{\tau_0})\right],
%\end{equation*}
\begin{multline} \label{eq:lastJ}
J^\nu(x)=\ME\left[1_{\{T_1<\tau_0\}}e^{-\alpha
T_1}\{K(X_{T_1-},X_{T_1})-g(X_{T_1}^0)+J^\nu(X_{T_1})\}\right]\\
+\ME\left[1_{\{T_1>\tau_0\}}e^{-\alpha\tau_0}\{P-g(X_{\tau_0})\}\right]+g(x).
\end{multline}
Noting that $g(X_{T_1}^0)=g(X_{T_1-}$, adding and subtracting
$g(X_{T_1})$ and further defining
\begin{equation*} \label{eq:uandJ}
u(x)\triangleq J^\nu(x)-g(x),
\end{equation*}
(\ref{eq:lastJ}) finally becomes
%a.s.
\begin{multline} \label{eq:simple2}
u(x)=\ME\left[1_{\{T_1<\tau_0\}}e^{-\alpha
T_1}\{K(X_{T_1-},X_{T_1})+u(X_{T_1})-g(X_{T_1-})+g(X_{T_1})\}\right]\\
+\ME\left[1_{\{T_1>\tau_0\}}e^{-\alpha\tau_0}\{P-g(X_{\tau_0})\}\right],
\end{multline}
and we consider the maximization of this $u(\cdot)$ function and
add back $g(x)$ since $\sup u(x)=\sup\{J^\nu(x)-g(x)\}=\sup
J^\nu(x)-g(x)$.  Note that this simplification leading to (\ref{eq:simple2}) does not depend on the $F$-concavity assumption.\\

Since we have confirmed that optimal intervention times are exit
times of the process from an interval, let us use a simpler notation
that $X_{T_i-}=b_i$ and $X_{T_i}=a_i$ for all $i$.  We can denote
$T_i-=\tau_b\triangleq \inf\{t>0; X_t\geq b_i\}$. By observing
(\ref{eq:simple2}),
\begin{align}\label{eq:a+b}
u(b)&=u(X_{T_1-})=K(X_{T_1-},X_{T_1})+g(X_{T_1})-g(X_{T_1-})+u(X_{T_1})\nonumber\\
&=K(b,a)+g(a)-g(b)+u(a)=\bar{K}(b, a)+u(a) \\
u(0)&=u(X_{\tau_0})=P-g(X_{\tau_0})=P-g(0)\nonumber
\end{align}
we have
\begin{eqnarray}\label{eq:twoside}
u(x)&=
\begin{cases}
u_0(x), &x\in[0, b)\\
\bar{K}(x, a)+u_0(a), &x\in[b, \infty).
\end{cases}
\end{eqnarray}
where
\begin{equation*}
u_0(x)\triangleq\ME[1_{\{\tau_b<\tau_0\}}e^{-\alpha \tau_b}u(b)]
+\ME[1_{\{\tau_b>\tau_0\}}e^{-\alpha\tau_0}u(0)].
\end{equation*}
\\
\noindent The second equation of (\ref{eq:twoside}) is obtained from
(\ref{eq:simple2}) by noticing that, on $x\in [b, \infty)$,
$\p^x(T_1<\tau_0)=1$.  Indeed, in this case, we immediately jump to
$a$, so that $X_{T_1-}=x$ and $X_{T_1}=a$. Since $a\in(0,b)$,
$u(a)=u_0(a)$.  Now let us note that we have the following
representations in (\ref{eq:simple2})
\begin{equation*}
\ME[e^{-\alpha\tau_r}1_{\{\tau_r<\tau_l\}}]=\frac{\psi(l)\varphi(x)-\psi(x)\varphi(l)}
{\psi(l)\varphi(r)-\psi(r)\varphi(l)}, \hspace{0.5cm} x\in[l,r]
\end{equation*}
where $\tau_l\triangleq\inf\{t>0; X_t=l\}$ and
$\tau_r\triangleq\inf\{t>0; X_t=r\}$ and $\varphi(\cdot)$ and
$\psi(\cdot)$ defined in the previous section.  Finally, with
$F(\cdot)$ being defined as in (\ref{eq:F}), we have a
characterization of $u(x)$,
\begin{equation}\label{eq:uchar}
\frac{u(x)}{\varphi(x)}=\frac{u(b)(F(x)-F(0))}{\varphi(b)(F(b)-F(0))}+\frac{u(0)(F(b)-F(x))}{\varphi(0)(F(b)-F(0))},
\quad x\in[0, b].
\end{equation}
\\
Define $W\triangleq\frac{u}{\varphi}\circ F^{-1}$, this becomes,
for any $a$ and $b$,
\begin{equation}\label{eq:W}
   W(F(x))=W(F(b))\frac{F(x)-F(0)}{F(b)-F(0)}+W(F(0))\frac{F(b)-F(x)}{F(b)-F(0)},
   \quad x\in[0, b].
\end{equation}
This represents a \emph{linear function} that passes a fixed point,
$A\triangleq(F(0), W(F(0))$.

To discuss how to find the optimal pair $(a^*, b^*)$, we write
$u(x)$ as $u_{a,b}(x)$ to emphasize the dependence on $a, b$, then
on $x\in[0, b]$,
\begin{align}
%\tilde{v}(x)&\triangleq
\sup_{a\in\mathbb{R_+}b\in\R_+}u_{a,b}(x)=\sup_{a\in\mathbb{R+}}\sup_{b\in\R_+}\{\ME[1_{\{\tau_b<\tau_0\}}e^{-\alpha
\tau_b}(\bar{K}(b,a)+u_{a,
b}(a))]+\ME[1_{\{\tau_b>\tau_0\}}e^{-\alpha\tau_0}u_{a,
b}(0)]\}.\label{eq:twostage}
\end{align}
This can be considered as a two-stage optimization problem. First,
let $a$ be fixed. For each $a$, the inner maximization of
(\ref{eq:twostage}) becomes
\begin{equation}\label{eq:stage1}
V_a(x)\triangleq\sup_{\tau_b\in\mathcal{S}}\{\ME[1_{\{\tau_b<\tau_0\}}e^{-\alpha
\tau_b}(\bar{K}(b,a)+V_a(a))]+\ME[1_{\{\tau_b>\tau_0\}}e^{-\alpha\tau_0}(P-g(0))]\}
\end{equation}
and, among $a's$, choose an optimal $a$ in the sense,
$\tilde{v}(x)\triangleq \sup_{a}V_a(x)$ for any $x$. It should be
pointed out that $V_a(x)$ may take negative values if $P-g(0)$ does.
Now, we discuss a solution method of the first stage optimization
(\ref{eq:stage1}).  For this purpose, we need a lemma:

%------------------
% Lemma 3-2
%=================
\begin{lemma}\label{lem:0}
If we define
\begin{equation*}
G(x,
\gamma)\triangleq\sup_{\tau\in\mathcal{S}}\ME[e^{-\alpha\tau}(h(X_\tau^0)+\gamma)],
\quad x\in\mathbb{R}, \gamma\in\mathbb{R}
\end{equation*}
for some Borel function $h: \mathbb{R}\rightarrow\mathbb{R}$ and
with condition (\ref{eq:fcon}), then, for $\gamma_1>\gamma_2\geq 0$,
\begin{equation} \label{eq:contraction}
G(x, \gamma_1)-G(x, \gamma_2)\leq \gamma_1-\gamma_2,
\end{equation}
for any $x$.
\end{lemma}
\begin{proof}
The left hand side of (\ref{eq:contraction}) is well-defined due
to (\ref{eq:fcon}). It is clear that $G(x, \gamma)$ is convex in
$\gamma$ for any $x$. Then
$D_\gamma^+G(x,\gamma_0)\triangleq\lim_{\gamma_\downarrow\gamma_0}\frac{G(x,\gamma_0)-G(x,
\gamma)}{\gamma_0 -\gamma}$ exists at every
$\gamma_0\in\mathbb{R}$, and
\begin{equation}\label{eq:D+}
\frac{G(x, \gamma_1)-G(x, \gamma_2)}{\gamma_1-\gamma_2}\leq
D^+_{\gamma}G(x, \gamma_1).
\end{equation}
Consider the bound of $G(x,\gamma)$ for $x$ fixed:
\begin{align*}
G(x,\gamma)\leq \sup_{\tau\in\mathcal{S}}\ME
e^{-\alpha\tau}|h(X_\tau)|+|\gamma|\sup_{\tau\in\mathcal{S}}\ME
e^{-\alpha\tau}.
\end{align*}

%First, consider a positive $\gamma>0$.
The first term on the right
hand side is constant in $\gamma$ and the second term is linear in
$\gamma$ and the $\ME[e^{-\alpha\tau}]\leq 1$ for any
$\tau\in\mathcal{S}$. Due to the convexity of $G(x, \gamma)$ in
$\gamma$, for the above inequality to hold, $D^+_\gamma G(x,
\gamma)\leq 1$ for all $\gamma\in\mathbb{R}$.  On account of
(\ref{eq:D+}), we have (\ref{eq:contraction}).
%Finally, consider a negative $\gamma<0$.  Since $G(x, \gamma)$ is
%a monotone increasing function of $\gamma$ once we fix
%$x\in\mathbb{R}$, the condition $D^+_\gamma G(x, \gamma)\leq 1$ is
%also necessary to hold in the range of $\gamma<0$.  The rest of the
%proof is the same as in a positive $\gamma$.
\end{proof}

%-----------------------------
% Lemma 3-3
%=============================
Coming back to (\ref{eq:stage1}), we need some care because the
value function $V_a(x)$ contains its value at $a$, $V_a(a)$ in the
definitive equation. Let us consider a family of optimal stopping
problem parameterized by $\gamma \in\mathbb{R}$.
\begin{align}\label{eq:gamma}
V^\gamma_a(x)&
\triangleq\sup_{\tau\in\mathcal{S}}\left\{\ME[1_{\{\tau<\tau_0\}}e^{-\alpha
\tau}(\bar{K}(X_\tau,a)+\gamma)]\nonumber
+\ME[1_{\{\tau>\tau_0\}}e^{-\alpha\tau_0}(P-g(0))]\right\} \nonumber \\
&=\sup_{\tau\in\mathcal{S}}\ME[e^{-\alpha \tau}r^\gamma(X_\tau, a)]
%&=\sup_{b}\ME[1_{\{\tau_b<\tau_0\}}e^{-\alpha
%\tau_b}(K(b,a)+g(a)-g(b)+\gamma)]+\ME[1_{\{\tau_b>\tau_0\}}e^{-\alpha\tau_0}(P-g(0))]
\end{align}
where
\begin{eqnarray} \nonumber
    r^\gamma(x, a) &=&  \begin{cases}
                 P-g(0), &x=0, \\
                %K(a, a), &0<x<a\\
                \bar{K}(x,a)+\gamma, &x>0.
     \end{cases}
\end{eqnarray}
Obviously, this parameterized problem can be solved by using
Proposition \ref{prop:2} to \ref{prop:4}.  Now we link this
parameterized optimal stopping problem to (\ref{eq:stage1}).
\begin{lemma}\label{lem:1}
For $a>0$ given, if there exists a solution to (\ref{eq:gamma}),
then there always exists unique $\gamma$ such that
$\gamma=V^\gamma_a(a)$ holds, provided that (\ref{eq:Kneg}) holds.
\end{lemma}

%-----------------------------
% Proof of Lemma
%=============================

\begin{proof}  Without loss of generality, we need only to consider
the case where
\begin{equation}\label{eq:positive-condition}
\sup_{x\in\mathbb{R}}\bar{K}(x, a)>0
\end{equation}
for some $a>0$.  Indeed, suppose that there is no such  $a$ and let
us consider a sequence of optimal stopping scheme.  In each
iteration, the value function for the optimal stopping problem takes
negative values, so that $\phi_n(\cdot)<0$ for all $n$.  Then in the
next iteration, $\bar{K}(x, y)$ function will be shifted downwards,
leading to $\phi_{n+1}(\cdot)<0$. Hence the ``no interventions"
strategy is trivially optimal.

In (\ref{eq:gamma}), since $\gamma$ is some constant parameter, we
benefit from Proposition \ref{prop:2} and claim that $V^\gamma_a(x)$
is characterized as the smallest $F$-concave majorant of
$r^\gamma(\cdot, a)$ that passes $\left(F(0),
\frac{P-g(0)}{\varphi(0)}\right)$. In terms of the notation of
Proposition \ref{prop:4}, if we define $W^\gamma_a(\cdot)$ such that
\begin{equation*}
V^\gamma_a(x)=\varphi(x)W^\gamma_a(F(x)),
\end{equation*}
then $W^\gamma_a(\cdot)$ passes through the fixed point $A=(F(0),
W^\gamma_a(F(0))$ and is the smallest concave majorant of
$H^\gamma(\cdot, a)\triangleq
\frac{r^\gamma(F^{-1}(\cdot),a)}{\varphi(F^{-1}(\cdot))}
=\frac{\bar{K}(F^{-1}(\cdot),
a)}{\varphi(F^{-1}(\cdot))}+\frac{\gamma}{\varphi(F^{-1}(\cdot))}$.

Now fix $a$. Our approach here is by starting with $\gamma=0$, we
move $\gamma$ and evaluate $V^\gamma_a(a)$ and try to find $\gamma$
such that $\gamma=V^\gamma_a(a)$. Due to
(\ref{eq:positive-condition}), we have $W^0_a(F(a))>0$.  By the
monotonicity of $F$, it is equivalent to saying that
$V^0_a(a)>0=\gamma$. As $\gamma$ increases, $W^\gamma_a(F(a))$
increases monotonically by the right hand side of (\ref{eq:gamma}).
Lemma \ref{lem:0} implies that for $\gamma_1>\gamma_2\geq 0$,
\begin{equation}\label{eq:contraction2}
V_a^{\gamma_1}(x)-V_a^{\gamma_2}(x)\leq \gamma_1-\gamma_2
\end{equation}
for any $x\in\mathbb{R_+}$. Note that $W_a^\gamma(F(a))\geq
H^\gamma(F(a), a)$.  However, since $V$ has less than the linear
growth in $\gamma$ as demonstrated by (\ref{eq:contraction2}), there
is a certain $\gamma^{'}$ large enough such that
$W^\gamma_a(F(a))=H^\gamma(F(a), a)$ for $\gamma\geq \gamma^{'}$.
This implies
\begin{align*}
\varphi(a)W^{\gamma^{'}}_a(F(a))&=\varphi(a)H^\gamma(F(a), a)\\
\Leftrightarrow V_a^{\gamma'}(a)&=\bar{K}(a, a)+\gamma'<\gamma'
\end{align*}
where the inequality is due to the assumption (\ref{eq:Kneg}). For
this $\gamma^{'}$, we have $V^{\gamma^{'}}_a(a)<\gamma^{'}$.

The monotonicity and continuity of $W^\gamma(F(a))$ (due to the
convexity of $V_a^\gamma(\cdot)$) with respect to $\gamma$, together
with (\ref{eq:contraction2}), implies that, for any $a$, there
exists one and only one $\gamma$ such that $V^\gamma_a(a)=\gamma$.
\end{proof}

%------------------------
% Methodology
%========================
\subsection{Methodology to find $v(x)$ and $(a^*, b^*)$}
Using (\ref{eq:uchar}), namely the
characterization of $u_{a,b}$, we describe an optimization procedure
based on Proposition \ref{prop:3} and \ref{prop:4}.
\begin{enumerate}
\item Fix $a$.  Consider the function
\begin{equation} \label{eq:myR}
R(\cdot, a)\triangleq\frac{\bar{K}(F^{-1}(\cdot),
a)}{\varphi(F^{-1}(\cdot))}
\end{equation}
%where
%\begin{eqnarray} \label{eq:myra}
%    r (x, a) &\triangleq&  \begin{cases}
%                P-g(0), &x=0, \\
%                %K(a, a), &0<x<a\\
%                 \bar{K}(x, a), &x>0.
%     \end{cases}
%\end{eqnarray}
Define $W_a(\cdot)$ such that $V_a(x)=\varphi(x)W_a(F(x))$ and by
the characterization (\ref{eq:u}), it is a straight line with a
slope, say $\beta(a)$ and passes through
$(F(0),W_a(F(0))=\left(F(0), \frac{P-g(0)}{\varphi(0)}\right)$. We
can write the linear majorant, in general,
\begin{equation} \label{eq:W2}
W_a(y)=\beta(a)y+d.
\end{equation}
%where in the absorbing boundary case, $d=-\beta F(0)+ W_a(F(0))$.
\item  \underline{First stage optimization:}  For each slope
$\beta(a)$, we can calculate the value of $W_a(F(a))$, but we have
to find the $W_a(\cdot)$ function such that, at some point
$F(b(a))$, we have
\begin{equation*}
W_a(F(b))=R(F(b), a)+W_a(F(a))\frac{\varphi(a)}{\varphi(b)}.
\end{equation*}
where we write $b(a)\equiv b$ for notational simplicity.  This
requirement is equivalent to finding $\gamma$ in (\ref{eq:gamma}) in
Lemma \ref{lem:1} such that
\begin{equation*}
\frac{\gamma}{\varphi(a)}=W_a^\gamma(F(a)).
\end{equation*}
 By Proposition \ref{prop:exitchar}, $\mathrm{C}_a\triangleq(0,
b(a))$ is the continuation region.  If $R$ is a differentiable
function with respect to the first argument, we can find the optimal
point $b(a)$ analytically. In effect, it is to find a point $b(a)$
such that the linear majorant and the shifted function $R(F(x),
a)+W_a(F(a))\frac{\varphi(a)}{\varphi(x)}$ have a tangency point.
This is equivalent to calculating the maximum slope that majorizes
$R(F(x), a)$ after it is shifted.  Explicitly, we solve
\begin{equation}\label{eq:simplesystem}
\left(\frac{\bar{K}(b,
a)}{\varphi(b)}\right)'-\frac{\varphi'(b)\varphi(a)}{\varphi(b)^2}d=\beta(a)\left(F'(b)+\frac{\varphi'(b)\varphi(a)}{\varphi(b)^2
}F(a)\right)
\end{equation}
for $b(a)$ where $\beta(a)$ is
\begin{equation}\label{eq:beta}
\beta(a)=\frac{\varphi(b)R(F(b),
a)-d(\varphi(b)-\varphi(a))}{F(b)\varphi(b)-F(a)\varphi(a)}.
\end{equation}
For the absorbing boundary case, these equations can be easily
modified.  Let us denote $D\triangleq
W_a(F(0))=(P-g(0))/\varphi(0)$.  Then (\ref{eq:simplesystem}) and
(\ref{eq:beta}) become
\begin{equation}\label{eq:simplesystem-absorbing}
\left(\frac{\bar{K}(b,
a)}{\varphi(b)}\right)'-\frac{\varphi'(b)\varphi(a)}{\varphi(b)^2}D
=\beta(a)\left(F'(b)+\frac{\varphi'(b)\varphi(a)}{\varphi(b)^2}(F(a)-F(0))\right)
\end{equation}
and
\begin{equation}
\beta(a)=\frac{\varphi(b)R(F(b),
a)+D\varphi(a)}{(F(b)-F(0))\varphi(b)-(F(a)-F(0))\varphi(a)},
\end{equation}
respectively.
\item   \underline{Second stage optimization:}  Now, let $a$ vary
and choose, among $\beta(a)$, find $\beta^*$, if exists, and also
the corresponding $b(a^*)$ and $a^*$.  Due to the characterization
of the value function, these $a^*$ and $b^*\triangleq b(a^*)$ must
be the solution to (\ref{eq:impulsevalue}).  Suppose that
$\bar{K}(x, a)$ is a decreasing function of $a$. As $a$ becomes
closer to $0$, the quantity $\bar{K}(x, a)$ becomes larger, while
$W_a(F(a))$ smaller.  Hence we can expect the existence of $a^*$
that maximizes the slope $\beta$.

\end{enumerate}

\begin{remark}\normalfont
With respect to the third point of the proposed method above, we
should check if there exists a concave majorant as $a\downarrow 0$.
%Again, for the notational convenience, we assume that $F(x)=x$. Then
Namely, we consider whether
\begin{equation*}
\lim_{x\downarrow 0}(K(x, a)-(g(x)-g(a))+u(a))<P-g(0)
\end{equation*}
holds in the neighborhood of $a=0$.  Since $\lim_{x\downarrow
0}g(x)=\lim_{a\downarrow 0}g(a)$ and $\lim_{a\rightarrow
0}u(a)=u(0)=P-g(0)$ by the continuity of $u$, the last inequality
holds due to (\ref{eq:Kneg}).\myBox \end{remark}

%---------------------------
% Section 3-3 General Case
%=========================

\subsection{Characterization of the Intervention Times and the Value Function:
General Case} \label{sec:gen} Let us move on to a general case where
the mapping $x\rightarrow\frac{\bar{K}}{\varphi}(x)$
$:\mathbb{R_+}\rightarrow\mathbb{R_+}$ is not necessarily
$F$-concave. First, we extend Proposition \ref{prop:exitchar} to
characterize optimal intervention times.
%------------------------
% Proposition 3-3
%========================

\begin{proposition}\label{prop:genK}
The value function $v(x)$ for (\ref{eq:impulsevalue}) is given by
the smallest solution majorizing $g$ of $v-g=\mathcal{L}(v-g)$ and
optimal intervention times $T_i^*$ are given by exit times from an
interval
%\begin{equation*}
%T_i^*=\inf\{t>T_{i-1}^*; X_t\notin(0, b_i^*), \quad i=1,2,....\}
%\end{equation*}
if and only if, for all $y\in\mathbb{R_+}$,
\begin{equation}\label{eq:iff}
x\rightarrow\bar{K}(x, y) \quad \text{is continuous and} \quad
q\triangleq\limsup_{x\rightarrow\infty}D^-\left(\frac{\bar{K}}{\varphi}\circ
F^{-1}\right)(x) \text{ is finite.}.
\end{equation}
where $D^-f(x_0)\triangleq\limsup_{x\uparrow
x_0}\frac{f(x)-f(x_0)}{x-x_0}$.
\end{proposition}
\begin{proof}
For any given $a\in\R_+$, if we can find the smallest linear
majorant of $\frac{\bar{K}(F^{-1}(\cdot),
a)+\gamma}{\varphi(F^{-1}(\cdot))}$ for an arbitrary
$\gamma\in\R_+$, we can find $\gamma=\varphi(a)W_a(F(a))$ by Lemma
\ref{lem:1}. Due to the constancy of $\gamma$, it suffices to show
that condition (\ref{eq:iff}) is necessary and sufficient for the
existence of concave majorant of $\frac{\bar{K}}{\varphi}\circ
F^{-1}$ on $F(\mathcal{I})$. The sufficiency is immediate. For the
necessity, we assume that $q=+\infty$.  We can take a sequence of
points $\{x^k\}\subset\mathbb{R}$ such that $x^k\rightarrow\infty$
and $D^-\left(\frac{\bar{K}}{\varphi}\circ
F^{-1}\right)(x^k)\rightarrow\infty$ as $k\rightarrow\infty$.  If
necessary, by taking a subsequence, we can make this sequence
$\{x^k\}$ monotone.  Consider the smallest concave majorant of
$\frac{\bar{K}}{\varphi}\circ F^{-1}$ on $[F(0),F(x^k)]$.  Call it
$v^k(x)$. It is clear that $v^k(x)$ is monotone increasing in $k$
for all $x\in[F(0),F(x^k)]$. As $k\rightarrow\infty$,
$x^k\rightarrow\infty$ and $v(x)\geq v^k(x)$.  We thus have
$v(x)=\lim_{k\rightarrow \infty}v^k(x)=\infty$ for all
$x\in\mathbb{R_+}$.  There is no optimal intervention policy.
\end{proof}
Suppose that the $F$-concavity of the reward function is violated,
so that the intervention point may be multiple.  Let us consider a
strategy that we have two intervention points, $b_1$ and $b_2$ being
arbitrarily chosen such that $0<b_1<b_2$. We want to characterize
function $J^\nu(x)$ as in (\ref{eq:J}) again.
%\begin{equation}\label{twobarrierprob}
%J^\nu(x)=\ME\left[\int_0^{\tau_0}e^{-\alpha
%s}f(X_s)ds+e^{-\alpha\tau_0}P+\sum_{T_i<\tau_0}e^{-\alpha
%T_i}K(X_{T_i-},X_{T_i})\right].
%\end{equation}
Recall that there are no controls in a way that the process is
pulled up to avoid ruin.  In other words, $\p^x[\tau_0<\infty]=1$.
Assume, for the moment, that we always apply control at these
boundaries $b_1$ and $b_2$ and then, once applied, the process moves
to $a_1<b_1$ and $a_2<b_2$, respectively.

If we start with a point $x\in[0, b_1]$, the problem is equivalent
to the case we considered already, since the process cannot go
beyond the level $b_1$. Hence following (\ref{eq:twoside}), we
have for $x\in[0, b_1]$
\begin{align*}
J^\nu_1(x)&\triangleq\ME\left[\int_0^{\tau_0}e^{-\alpha
s}f(X_s)ds+e^{-\alpha\tau_0}P+\sum_{T_i<\tau_0}e^{-\alpha
T_i}K(X_{T_i-},X_{T_i})\right]
\end{align*}
and
\begin{align*}
u_1(x)&=\ME[1_{\{\tau_{b_1}<\tau_0\}}e^{-\alpha
\tau_{b_1}}u_1(b_1)]
+\ME[1_{\{\tau_{b_1}>\tau_0\}}e^{-\alpha\tau_0}u_1(0)], &x\in[0,
b_1]
\end{align*}
by defining $u_1(x)\triangleq J^\nu_1(x)-g(x)$.  If we start with
a point $x\in[b_1, b_2]$, there are two strategies available:
\begin{itemize}
    \item [(A)] Let $X_t$ move along.  (It either hits $b_1$ or $b_2$
    first.)
    \item [(B)] Apply the control immediately $(t=0)$ by moving the process from $x$ to $a_1$ (the post-control point that corresponds to $b_1$)
    and let the process start at $a_1$. (Recall that we do not let
    $X$ enter into $(b_1, \infty)$ after moving to $a_1$.)
\end{itemize}
Consider strategy (A) first. Let us define
\begin{equation*}
J^\nu_2(x)\triangleq\ME\left[\int_0^{\tau_{b_1}}e^{-\alpha
s}f(X_s)ds+\sum_{T_i<\tau_{b_1}}e^{-\alpha
T_i}K(X_{T_i-},X_{T_i})\right] \quad x\in[b_1, b_2].
\end{equation*}
Using the strong Markov property, we can reduce $J^\nu_2$ to a
simpler form. For any $(a_1, b_1)$ and $(a_2, b_2)$, we have
\begin{multline*}
   J_2^\nu(x)=\ME[1_{\{\tau_{b_1}<\tau_{b_2}\}}e^{-\alpha
\tau_{b_1}}K(X_{\tau_{b_1-}},
   X_{\tau_{b_1}})-g(X_{\tau_{b_1}})+J_1^\nu(X_{\tau_{b_1}})+g(X_{\tau_{b_1}})-g(X_{\tau_{b_1-}})\\
   +\ME[1_{\{\tau_{b_1}>\tau_{b_2}\}}e^{-\alpha
\tau_{b_2}}K(X_{\tau_{b_2-}},
   X_{\tau_{b_2}})-g(X_{\tau_{b_2}})+J_2^\nu(X_{\tau_{b_2}})+g(X_{\tau_{b_2}})-g(X_{\tau_{b_2-}})]+g(x).
\end{multline*}
We shall use $u_1(x)=J_1^\nu(x)-g(x)$ in the first term. Now let
us define $u_2(x)\triangleq J_2^\nu(x)-g(x)$. Then the last
equation becomes
\begin{align}\label{eq:uu}
    u_2(x)&=\ME[1_{\{\tau_{b_1}<\tau_{b_2}\}}e^{-\alpha
\tau_{b_1}}K(X_{\tau_{b_1-}},
  X_{\tau_{b_1}})+u_1(X_{\tau_{b_1}})+g(X_{\tau_{b_1}})-g(X_{\tau_{b_1-}})\nonumber\\
&\hspace{2cm}+\ME[1_{\{\tau_{b_1}>\tau_{b_2}\}}e^{-\alpha
\tau_{b_2}}K(X_{\tau_{b_2-}},
   X_{\tau_{b_2}})+u_2(X_{\tau_{b_2}})+g(X_{\tau_{b_2}})-g(X_{\tau_{b_2-}})]\nonumber\\
&=\ME[1_{\{\tau_{b_1}<\tau_{b_2}\}}e^{-\alpha
\tau_{b_1}}(\bar{K}(b_1,a_1)+u_1(a_1))]+\ME[1_{\{\tau_{b_1}>\tau_{b_2}\}}e^{-\alpha
\tau_{b_2}}(\bar{K}(b_2, a_2)+u_2(a_2))]
\end{align}
on $x\in[b_1, b_2]$.  By identifying $\bar{K}(b_2,
a_2)+u_2(a_2)=u_2(b_2)$ and
$u_2(b_1)=\bar{K}(b_1,a_1)+u_1(a_1)=u_1(b_1)$ that shows $u_1(x)$
and $u_2(x)$ are connected at $x=b_1$. Thus,
\begin{equation} \label{twobarrier}
u_2(x)=\frac{\varphi(x)}{\varphi(b_1)}\frac{F(b_2)-F(x)}{F(b_2)-F(b_1)}u_2(b_1)
+\frac{\varphi(x)}{\varphi(b_2)}\frac{F(x)-F(b_1)}{F(b_2)-F(b_1)}u_2(b_2),
\quad x\in[b_1, b_2].
\end{equation}
To summarize this result, if we define
$W_i(\cdot)\triangleq\frac{u_i}{\varphi}\circ F^{-1}(\cdot)$ for
$i=1,2$ on $F(\mathcal{I})$, this is again a linear function for
each $i$. Hence by defining
\begin{eqnarray} \nonumber
    W_A(F(x))\triangleq  \begin{cases}
                 W_1(F(x))=W_1(F(0))\frac{F(b_1)-F(x)}{F(b_1)-F(0)}+W_1(F(b_1))\frac{F(x)-F(0)}{F(b_1)-F(0)}, \hspace{0.4cm}x\in[0, b_1] \\
                W_2(F(x))=W_2(F(b_1))\frac{F(b_2)-F(x)}{F(b_2)-F(b_1)}+W_2(F(b_2))\frac{F(x)-F(b_1)}{F(b_2)-F(b_1)},
                 \hspace{0.2cm}x\in[b_1, b_2],
     \end{cases}
\end{eqnarray}
we have a piecewise linear function on $F(\mathcal{I})$.  Moreover,
since we can treat $b_1$ as an absorbing boundary, we have
$a_1<b_1<a_2<b_2$.

Next consider strategy (B), whose value function is
\begin{eqnarray}\label{eq:W-B}
W_B(F(x))\triangleq \begin{cases} W_1(F(x)), &0\leq x \leq b_1\\
\overline{W_1}(F(x))\triangleq
\frac{\varphi(a_1)}{\varphi(x)}W_1(F(a_1))+R(F(x), a_1), &b_1<x.
\end{cases}
\end{eqnarray}

%=========================
% Lemma 3-3
%=========================
\begin{lemma} \label{lem:2}
(A) is better than (B) only if
\begin{align*}
\beta_1\triangleq\frac{W(F(b_1))-W(F(0))}{F(b_1)-F(0)} <
\frac{W(F(b_2))-W(F(b_1))}{F(b_2)-F(b_1)}\triangleq \beta_2.
\end{align*}
\end{lemma}
\begin{proof}
Since the value function of strategy (B) is (\ref{eq:W-B}), choosing
(A) over (B) is equivalent to
\begin{equation*}
\overline{W_1}(F(x))<W_2(F(x)) \quad \text{on}\hspace{0.2cm}x>b_1.
\end{equation*}
If $W_1(F(x))$ majorizes $\overline{W_1}(F(x))$ on $x\in[0,
\infty)$, then this problem reduces to $F$-concavity case
discussed in the previous subsection.  Hence we consider the case
where there exists some $x\in[b_1, \infty)$ such that
\begin{equation*}
W_1(F(x))<\overline{W_1}(F(x)).
\end{equation*}  Now suppose that we have $\beta_1\leq \beta_2$. Then it is clear that
we cannot have $W_2(F(x))>\overline{W_1}(F(x))$ on $x\in[b_1,
\infty)$.
\end{proof} \noindent There are two cases to consider:
\begin{itemize}
\item [(1)] If $W_2(F(x))$ majorizes $\overline{W_1}(F(x))$ on
$x\in [b_1, \infty)$, then we adopt the point $b_2$ as an
intervention point. In this case, $\beta_2>\beta_1$ holds. However,
this implies that if we connect $A\triangleq(F(0), W_1(F(0))$ and
$C\triangleq(F(b_2), W_2(F(b_2))$, then this line segment $AC$ is
above the line segment connecting, piece by piece, points $A$,
$B\triangleq(F(b_1), W_1(F(b_1))$ and $C$. We can show that there
exists a point $b'\geq b_2$ such that its corresponding linear
majorant $W'(F(x))$ satisfies $W'(F(x))>W_1(F(x))$ on $x\in[0, b_1]$
and $W'(F(x))>W_2(F(x))$ on $[b_1, b_2]$. The proof of the existence
of a post-intervention point $a'$ corresponding to this point $b'$
follows in a similar manner to Lemma \ref{lem:1}.
\item[(2)] If $W_2(F(x))$ does not majorize
$\overline{W_1}(F(x))$, we can find another point $\bar{b}$, instead
of $b_1$, such that the linear (not piecewise linear) function
$W(F(x))$ corresponds to $\bar{b}$ majorizes $R(F(x),
\bar{a})+W(F(\bar{a}))\frac{\varphi(\bar{a})}{\varphi(x)}$ on $x\in
\mathbb{R}_+$ by Proposition \ref{prop:genK}.
\end{itemize}
In either case, the value function in the transformed space should
be a linear function that attains the largest slope among all the
possible linear majorant.  This argument holds true for any $b_1$
and $b_2$ with $b_1<b_2$.  We can continue this argument inductively
to the case of $n$ intervention points, $(b_1,...b_n)$. %For a local
%candidate $b_i^*$ to be chosen as a barrier point, the slope of
%$W_i(\cdot)$, denoted as $\beta_i$, must be greater than $\beta_1$,
%$\beta_2$,..., $\beta_i$. But this implies that if we connect
%$(F(0), W(F(0))$ and the $(F(b_i^*), W_i(F(b_i^*))$, that line
%segment is above the line segment connecting, piece by piece, all
%the points
%\begin{equation*}
%(F(0), W_1(F(0)), (F(b_1^*), W_1(F(b_1^*)),......,(F(b_i^*),
%W_i(F(b_i^*)).
%\end{equation*}
%Now, using the argument for the two points case $(b_1^*, b_2^*)$,
%it follows that there exists a point $x=b'$ such that the
% line segment $W'(F(x))$ connecting $(F(0),
%W_1(F(0))$ and $(F(b'), W'(F(b'))$ satisfies
%\begin{equation*}
%W'(F(x))>W(F(x)), \quad \text{for} \hspace{0.2cm} \forall
%x\in[0,\infty)
%\end{equation*}
%for all the piecewise linear functions $W(F(x))$ and that
%$W'(\cdot)$ is the value function (\ref{eq:impulsevalue}). The
%construction of $W'(\cdot)$ is given in the previous subsection.
%----------------------------------------------------------
We here summarize our argument up to this point as a main
proposition:
%-----------------------------
% Proposition
%=============================
\begin{proposition}
Suppose that (\ref{eq:iff}) holds and the optimal continuation
region is connected. The value function corresponding to
(\ref{eq:J}) of the impulse problem described in
(\ref{eq:Z})$\sim$ (\ref{eq:impulsevalue}) is written as
\begin{eqnarray}\label{eq:sol1}
v(x)&=& \begin{cases} v_0(x)\triangleq\varphi(x)W^*(F(x))+g(x), & 0\leq x\leq b^*\\
                        v_0(a^*)+K(x, a^*), &b^*\leq x.
        \end{cases}
\end{eqnarray}
where $W^*(\cdot)$ is the line segment connects $(F(0), W^*(F(0)))$
and $(F(b^*), W^*(F(b^*)))$ and satisfy the following:
\begin{enumerate}
    \item $W^*(F(\cdot))$ is the smallest linear majorant of $W^*(F(a^*))\frac{\varphi(a^*)}{\varphi(\cdot)}+R(F(\cdot),
    a^*)$
     and
    meets with $W^*(F(a^*))+R(F(\cdot), a^*)\frac{\varphi(a^*)}{\varphi(\cdot)}$ at point $F(b^*)$ and passes $(F(0),
\frac{P-g(0)}{\varphi(0)})$.  If $R$ is differentiable,  $(a^*,
    b^*)$ satisfy (\ref{eq:simplesystem}).
 \item The slope of $W^*(\cdot)$, denoted as $\beta^*$, is the
    largest slope among $\beta (a)$'s of all the possible linear majorants $W_a(\cdot)$.
\end{enumerate}
Moreover, if the mapping $x\rightarrow\frac{\bar{K}}{\varphi}(x)$
$:\mathbb{R}_+\rightarrow\mathbb{R}_+$ is $F$-concave, then the
optimal continuation region $(0, b^*)$ is uniquely determined.
\end{proposition}
Note that, at $x=0$,
\begin{equation*}
v(0)=\varphi(0)W^*(F(0))+g(0)=\varphi(0)\frac{P-g(0)}{\varphi(0)}+g(0)=P
\end{equation*}
as expected.

\begin{remark}\normalfont \label{rmk:pathology}
If the $F$-concavity of $\bar{K}$ is violated, there are two
possible cases (and combination of them) of multiple continuation
regions.
\begin{enumerate}
    \item For some $a_i^*$ with $i=1, 2,...$, we have the
    common $\beta^*$.  This is the case which we shall show in the
    next example.  In this case, the continuation region is $\mathrm{C}=\{(0, b_1^*), (b_1^*, b_2^*), (b_2^*, b_3^*) ...\}$
where $b_i^*$ corresponds to $a_i^*$ for each $i$, and the
    intervention region is $\Gamma= \{\{b_1^*\}, \{b_2^*\},
    \{b_3^*\}...\}$.  Each time the process hits one of the points $\{b_i^*\}$, the
control pulls the process back to the corresponding $a_i^*$.
    \item Another case is that, for the unique optimal $a^*$,
    there exists non-unique $b_1^*$ and $b_2^*$. In this case, the
    continuation region is $\mathrm{C}=\{(0, b_1^*), (b_1^*, b_2^*)\}$,
and the stopping region is $\Gamma=\{\{b_1^*\}, [b_2^*, \infty)\}$.
If the process hits $b_1^*$ or $b_2^*$, then the control pulls the
process back to $a^*$ in either situation.  It makes sense to
continue in the region $(b_1^*, b_2^*)$ because there is a positive
probability that one can extract $\bar{K}(b_2^*,
a^*)(>\bar{K}(b_1^*, a^*))$ within a finite time.
\end{enumerate}
\end{remark}

%----------------------------
%  Section 4  Free Boundary
%============================
\subsection{No absorbing boundary case} Next, we extend our argument to a
problem without the absorbing boundary. Hence the process can move
along in the state space in an infinite amount of time.  The problem
becomes
\begin{equation} \label{eq:fJ}
    J^\nu(x)=\ME\left[\int_0^{\infty}e^{-\alpha s}f(X_s)ds+\sum_{i=1}e^{-\alpha T_i}K(X_{T_i-},X_{T_i})\right]
\end{equation}
%\begin{equation} \label{eq:impulsevalue2}
%v(x)\triangleq\sup_{\nu\in\mathcal{V}} J^{\nu}(x)=J^{\nu^*}(x).
%\end{equation}
We can characterize intervention times as exit times from certain
boundary and simplify the performance measure (\ref{eq:fJ})
\begin{align*}
J^\nu(x)&=\ME\left[\int_0^\infty e^{-\alpha
s}f(X_s)ds+\sum_{i=1}e^{-\alpha
T_i}K(X_{T_i-},X_{T_i})\right]\nonumber\\
&=\ME[e^{-\alpha T_1}\{K(X_{T_1-},
X_{T_1})-g(X_{T_1-})+J^\nu(X_{T_1})\}]+g(x).
\end{align*}
The second equation is easily obtained in the same way as in the
previous section by noting $\p^x(T_1<\infty)=1$. The last term does
not depend on controls, so we define $u(x)\triangleq J^\nu(x)-g(x)$:
\begin{equation*} \label{eq:u3}
u(x)=\ME[e^{-\alpha T_1}\{K(X_{T_1-},
X_{T_1})-g(X_{T_1-})+g(X_{T_1})+u(X_{T_1})\}].
\end{equation*}
Again, we consider the $F$-concave case with the notation
$T_i-=\tau_b$ for all $i$ and we have
\begin{align*}
u(x)=\ME[e^{-\alpha \tau_b}(K(b,a)-g(b)+g(a)+u(a))]=\ME[e^{-\alpha
\tau_b}(\bar{K}(b,a)+u(a))].
\end{align*}
By defining $W=(u/\varphi)\circ F^{-1}$, we have
\begin{align*}\label{eq:W0}
W(F(x))&=W(F(c))\frac{F(b)-F(x)}{F(b)-F(c)}+W(F(b))\frac{F(x)-F(c)}{F(b)-F(c)},
\quad x\in(c, b].
\end{align*}
We should note that $F(c)\triangleq F(c+)=\psi(c+)/\varphi(c+)=0$
and
\begin{equation*}
W(F(c))=l_c\triangleq \limsup_{x\downarrow c}\frac{\bar{K}(x,
a)^+}{\varphi(x)}
\end{equation*}
for any $a\in (c,d]$.  For more detailed mathematical meaning of
this value $l_c$, we refer the reader to Dayanik and
Karatzas\cite{DK2003}.  We can effectively consider $(F(c), l_c)$ as
the absorbing boundary.

%--------------------------
% Section Example
%==========================
\section{Examples}
In this section, we work out some examples from financial
engineering problems.  For this purpose, we recall some useful
observations.  If $h(\cdot)$ is twice-differentiable at $x\in
\mathcal{I}$ and $y\triangleq F(x)$, then $H^{'}(y)=m(x)$ and
$H^{''}(y)=m^{'}(x)/F^{'}(x)$ with
\begin{equation}\label{eq:devH}
m(x)= \frac{1}{F^{'}(x)}\left(\frac{h}{\varphi}\right)^{'}(x),
\quad \text{and} \quad H^{''}(y) (\mathcal{A}-\alpha)h(x)\geq 0,
\quad y=F(x)
\end{equation}
with strict inequality if $H^{''}(y)\neq 0$.  These identities are
of practical use in identifying the concavities of $H(\cdot)$ when
it is hard to calculate its derivatives explicitly.
%---------------------------------
%  Example 1 Oksendal 1999
%=================================
\begin{example} \normalfont \label{ex:oks}
{\O}ksendal \cite{O1999} considers the following problem:
\begin{equation}\label{oksprob}
   J^\nu_o(x)=\ME\left[\int_0^\infty e^{-\alpha s}X_s^2 ds +
   \sum_{i}^{\infty}e^{-\alpha T_i}(c+\lambda\xi_i)\right]
\end{equation}
where $X^0_t=B_t$ is a standard Brownian motion and $c>0$ and
$\lambda\geq 0$ are constants.  The Brownian motion represents the
exchange rate of some currency and each impulse represents an
interventions taken by the central bank in order to keep the
exchange rate in a given target zone. Here we are only allowed to
give the system impulses $\zeta$ with values in $(0,+\infty)$. By
reducing a level from $b$ to $a$ (i.e., $b>a$) through
interventions, one can save continuously incurred cost (which is
high if the process is at a high level). The problem is to minimize
the expected total discounted cost $v_o(x)=\inf_{\nu} J_o^\nu(x)$.
We want to solve its sup version and change the sign afterwards
(i.e. $v_o(x)=-v(x)$):
\begin{equation*}
v(x)=\sup_{\nu}\ME\left[\int_0^\infty e^{-\alpha s}(-X_s^2) ds -
   \sum_{i}^{\infty}e^{-\alpha T_i}(c+\lambda\xi_i)\right].
\end{equation*}
The continuous cost rate $f(x)=-x^2$ and the intervention cost is
$K(x, y)=-c-\lambda(x-y)$ in our terminology. By solving the
equation $(\mathcal{A}-\alpha)v(x)=\frac{1}{2}v^{''}(x)-\alpha
v(x)=0$, we find $\psi(x)=e^{x\sqrt{2\alpha}}$ and
$\varphi(x)=e^{-x\sqrt{2\alpha}}$.  Hence
$F(x)=e^{2x\sqrt{2\alpha}}$ and $F^{-1}(x)=\frac{\log
x}{2\sqrt{2\alpha}}$.  Following our characterization of the value
function, we obtain
\begin{align*}
    J^\nu(x)&=\ME[e^{-\alpha T_1}\left\{K(X_{T_1-}, X_{T_1})-g(X^0_{T_1-})+J^\nu(X_{T_1})\right\}]+g(x)
\end{align*}
where $g(x)$ can be calculated by Fubini's theorem:
\begin{align*}
g(x)=-\ME\int_0^\infty e^{-\alpha
s}(x+B_s)^2ds=-\left(\frac{x^2}{\alpha}+\frac{1}{\alpha^2}\right).
\end{align*}
By defining $u(x)=J^\nu(x)-g(x)$, we have $u(x)=\ME[e^{-\alpha
\tau_b}\left\{K(b, a)-g(b)+g(a)+u(a)\right\}]$. Note that when
$b>a$, $g(a)-g(b)>0$ is the source of cost savings.

Let us fix $a>0$ and consider
$h(x)\triangleq-c-\lambda(x-a)+\frac{x^2-a^2}{\alpha}$ and
$H(y)\triangleq(h/\varphi)(F^{-1}(y)), y>0$.  By the first
equation in (\ref{eq:devH}), the sign of
$\left(\frac{h}{\varphi}\right)'(x)$ will lead us to conclude that
$H(F(x))$ is increasing from a certain point, say $x=p$ on $(p,
\infty)$, so is $H(F(x))$. Also, by direct calculation,
$H'(+\infty)=0$, from which we can assert that the value function
is finite by Proposition \ref{prop:genK}.  If we set
$p(x)\triangleq-x^2+a^2+\lambda\alpha(x-a)+\alpha c+1/\alpha$,
then $(\mathcal{A}-\alpha)h(x)=p(x)$ for every $x>0$. This
quadratic function $p(x)$ possibly has one or two positive roots.
 Let $k$ be the largest one.   Since $\lim_{x\rightarrow \infty}p(x)=-\infty$,
 by the second inequality in
(\ref{eq:devH}), $H(\cdot)$ is concave on $(F(k), +\infty)$. Hence
$H(y, a)$ is increasing and concave on $y\in (F(k), \infty)$.  %It
%can be shown that $H(F(p))<0$, $H(F(a))<0$ and $p<k$. From these
%facts,
Since the cost function in the transformed space is increasing and
concave from a certain point on, there is a linear majorant that
touches the cost function once and only once.  We can conclude that
for any $a>0$ and the parameter set, we have a connected
continuation region in the form of $(0, b^*)$.

For this fixed $a$, let us define $W_a(\cdot)$ such that
$V_a(x)=\varphi(x)W_a(F(x))$ and $r(x, a)=-c$ if $x<a$ and
$-c-\lambda(x-a)+\frac{x^2}{\alpha}-\frac{a^2}{\alpha}$ if $x\geq
a$. Then we have for any $a>0$,
\begin{equation}\label{}
    l_{-\infty}=\limsup_{x\downarrow-\infty}\frac{r(x, a)^+}{\varphi(x)}=0.
\end{equation}
Recall that the left boundary $-\infty$ is natural for a Brownian
motion.  Hence $W_a(y)$ that passes the origin of the transformed
space is the straight-line majorant of $R(\cdot,
a)+W_a(F(a))/\varphi(F^{-1}(\cdot))$ where $R(\cdot, a)$ is defined
in (\ref{eq:myR}):
\begin{eqnarray*}
R(y, a)= \begin{cases}

                    -c\sqrt{y}, &0 \leq y\leq F(a), \\
                       H(y, a)=\sqrt{y}\left(-c
                       -\frac{\lambda}{2\sqrt{2\alpha}}\log y+\lambda a
                       +\frac{(\log y)^2}{8\alpha^2}-\frac{a^2}{\alpha}\right),
                       &y>F(a).
        \end{cases}
\end{eqnarray*}
We can represent $W_a$ as $W(y)=\beta y$.  Since $R(x, a)$ is
differentiable with respect to $x$ on $x\geq a$, we can use
(\ref{eq:simplesystem}) to find $b(a)$ and corresponding
$\beta(a)$.  Then varying $a$, one can find the optimal $(a^*,
b^*, \beta^*)$.  Going back to the original space, on
$x\in(-\infty, b^*]$
\begin{equation*}
\tilde{v}(x)\triangleq \sup
u(x)=\varphi(x)W^*(F(x))=\varphi(x)(\beta^*)F(x)=\beta^*
e^{x\sqrt{2\alpha}}.
\end{equation*}  To get $v(x)=\sup_{\nu}J^\nu(x)$, we add back $g(x)$,
\begin{equation*}
v(x)=\tilde{v}(x)+g(x)=\beta^*e^{x\sqrt{2\alpha}}-\left(\frac{x^2}{\alpha}+\frac{1}{\alpha^2}\right).
\end{equation*}
Finally, flip the sign and obtain the optimal cost function
\begin{eqnarray*}
v_o(x)&=& \begin{cases}
                    \hat{v}_o(x)\triangleq\left(\frac{x^2}{\alpha}+\frac{1}{\alpha^2}\right)-\beta^*e^{x\sqrt{2\alpha}}, &0\leq x \leq b^*, \\
                    \hat{v}_o(a^*)+c+\lambda(x-a^*). &b^*\leq x.
        \end{cases}
\end{eqnarray*}
which coincides with the solution given by {\O}ksendal \cite{O1999}.
Figure 1 displays the solution with parameters $(c, \lambda,
\alpha)=(150, 50, 0.2)$.
\end{example}
\begin{figure}[h]%\label{fig:forex}
\begin{center}
\begin{minipage}{0.45\textwidth}
\centering \includegraphics[scale=0.75]{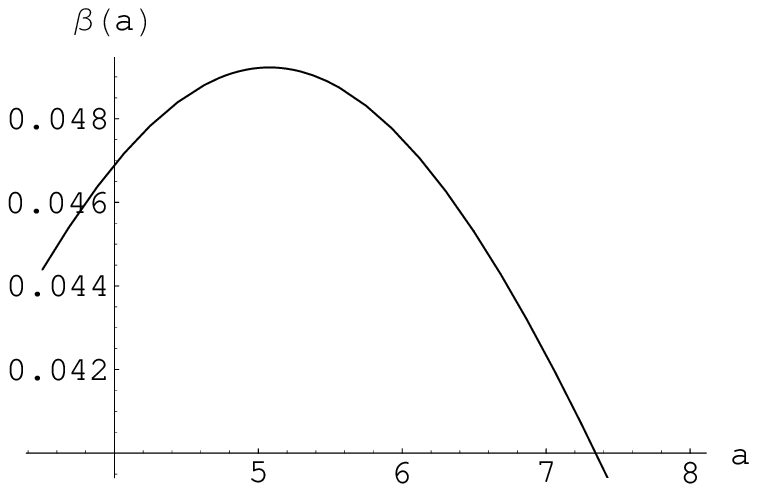} \\
(a)
\end{minipage}
\begin{minipage}{0.45\textwidth}
\centering \includegraphics[scale=0.75]{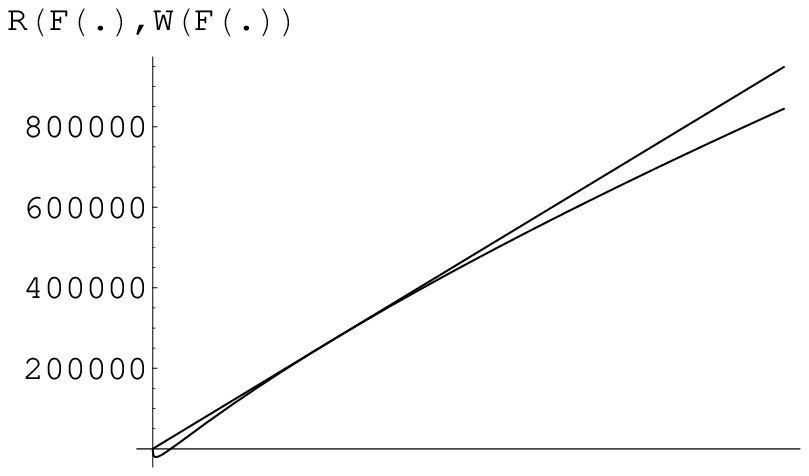} \\
(b)
\end{minipage}
\begin{minipage}{0.45\textwidth}
\centering \includegraphics[scale=0.75]{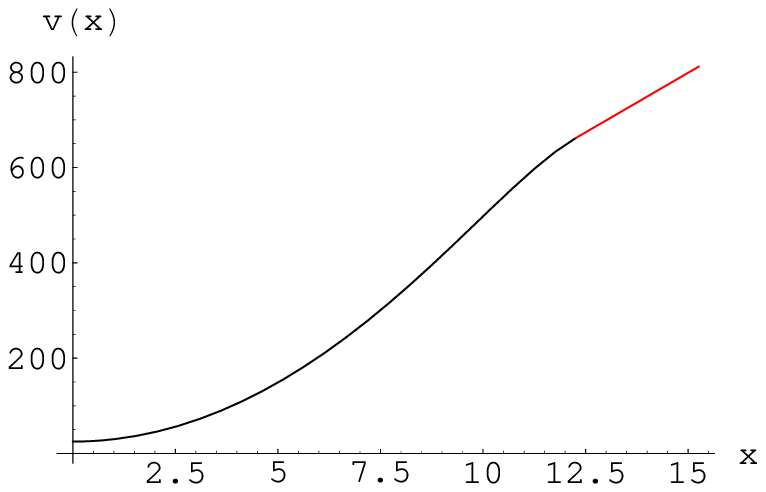} \\
(c)
\end{minipage}
\begin{minipage}{0.45\textwidth}
\centering \includegraphics[scale=0.75]{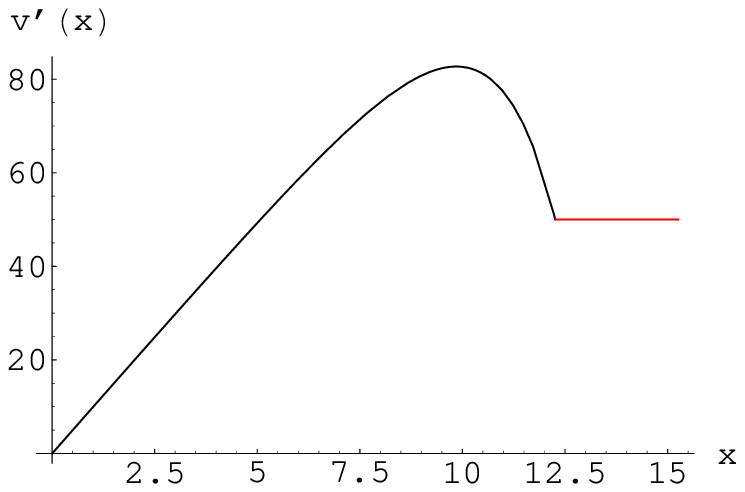} \\
(d)
\end{minipage}
\caption{\small(a) The plot of $\beta(a)$ against $a$, the former
being maximized at $a^*=5.077$ with $\beta^*=0.0492$. (b) The
functions $R(F(\cdot), a^*)$ shifted by the amount
$W_{a^*}(F(a^*))\frac{\varphi(a)}{\varphi(x)}$ (lower curve) and the
majorant $W_{a^*}(F(\cdot))$ (upper curve) corresponding to $a^*$,
giving us $b^*=12.261$. (c) The cost function $v_o(x)$. (d) The
derivative of $v_o(x)$, showing that the smooth-fit principle holds
at $b^*$.}
\end{center}
\end{figure}

\begin{figure}[h]\label{fig:2}
\begin{center}
\begin{minipage}{0.45\textwidth}
\centering \includegraphics[scale=0.75]{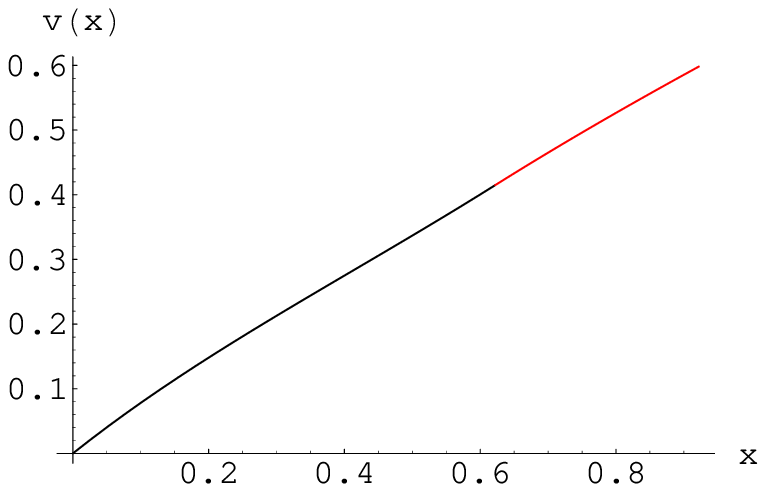} \\
\end{minipage}
\caption{\small The value function for Cadenillas et al.\cite{CSZ2003} problem.}
\end{center}
\end{figure}

%----------------------------
%  Example OU
%============================
\begin{example}\normalfont \label{ex:OU}
This example is a dividend payout problem where the underlying
process follows an Ornstein-Uhlenbeck process.  This problem was
originally studied by Cadenillas et al.\cite{CSZ2003} in an
ingenious way, but the existence of the finite value function and
the connectedness of the continuation region were left open.
Suppose that $X^0$ has the dynamics
\begin{equation*}
dX_t^0= \delta(m-X_t)dt +\sigma dW_t, \quad t\geq 0,
\end{equation*}
where $\delta>0$, $\sigma >0$ and $m\in\mathbb{R}$.  Only positive
impulse is allowed in this problem.  We consider the impulse control
problem,
\begin{equation*}
v(x)\triangleq\sup_{\nu\in
S}\ME\left[\sum_{T_i<\tau_0}^{\infty}e^{-\alpha
T_i}(-K+k\xi_i^\gamma)\right].
\end{equation*}
with some positive constant $K$, $k$ and the risk-aversion parameter
$\gamma\in(0, 1]$.  Since $\xi\in\mathbb{R_+}$, we have
\begin{equation*}
\bar{K}(x, y) = k(x-y)^\gamma-K, \quad x>y>0.
\end{equation*}
The functions $\psi(\cdot)$ and $\varphi(\cdot)$ are positive,
increasing and decreasing solutions of the differential equation
$(\mathcal{A}-\alpha)v(x)=(1/2)\sigma^2
v^{''}(x)+\delta(m-x)v^{'}(x)-\alpha v(x)=0$.  We denote, by
$\tilde{\psi}(\cdot)$ and $\tilde{\varphi}(\cdot)$, the functions of
the fundamental solutions for the auxiliary process
$Z_t\triangleq(X_t-m)/\sigma, t\geq 0$, which satisfies
$dZ_t=-\delta Z_t dt+dW_t$.  For every $x\in \mathbb{R}$,
\begin{equation*}
\tilde{\psi}(x)=e^{\delta
x^2/2}\mathcal{D}_{-\alpha/\delta}(-x\sqrt{2\delta})\quad\text{and}\quad
\tilde{\varphi}(x)=e^{\delta
x^2/2}\mathcal{D}_{-\alpha/\delta}(x\sqrt{2\delta}),
\end{equation*}
and $\psi(x)=\tilde{\psi}((x-m)/\sigma)$ and
$\varphi(x)=\tilde{\varphi}((x-m)/\sigma)$, where
$\mathcal{D}_\nu(\cdot)$ is the parabolic cylinder function; (see
Borodin and Salminen (2002, Appendices 1.24 and 2.9) and Carmona
and Dayanik (2003, Section 6.3)).  By using the relation
\begin{equation} \label{eq:Dft}
\mathcal{D}_\nu(z)=2^{-\nu/2}e^{-z^2/4}\mathcal{H}_\nu(z/\sqrt{2}),
\quad z\in\mathbb{R}
\end{equation}
in terms of the Hermite function $\mathcal{H}_\nu$ of degree $\nu$
and its integral representation
\begin{equation}\label{eq:Hermite}
\mathcal{H}_\nu(z)=\frac{1}{\Gamma(-\nu)}\int_0^\infty
e^{-t^2-2tz}t^{-\nu-1}dt, \quad \text{Re}(\nu)<0,
\end{equation}
(see for example, Lebedev(1972, pp284, 290)).  Let us consider the
function
\begin{equation*}
h(x)\triangleq kx^\gamma-K, \quad x>0, \gamma\in(0,1].
\end{equation*}
Since the function $h(\cdot)$ is increasing, the function
$H(y)=(h/\varphi)\circ F^{-1}(y), y\in(0, \infty)$ is also
increasing. Let us define the function
\begin{equation*}
p(x)\triangleq\frac{1}{2}\sigma^2
k\gamma(\gamma-1)x^{\gamma-2}+m\delta k\gamma
x^{\gamma-1}-k(\delta \gamma+\alpha)x^{\gamma}+\alpha K
\end{equation*}
which satisfies $(\mathcal{A}-\alpha)h(x)=p(x)$.  By using
(\ref{eq:devH}), $H^{''}(y)$ and $p(F^{-1}(y))$ have the same sign
at every $y$ where $h$ is twice-differentiable.  Hence we study the
(positive) roots of $p(x)=0$.  We have to divide two cases: (1)
$\gamma=1$ and (2) $\gamma<1$.  In either case, it can be shown that
$H^{'}(+\infty)=0$ by using (\ref{eq:Dft}) and (\ref{eq:Hermite})
and the identity $\mathcal{H}'_\nu(z)=2\nu\mathcal{H}_{\nu-1}(z),
z\in\mathbb{R}$.  Therefore, by Proposition \ref{prop:genK}, the
finiteness of the value function is proved.
\begin{itemize}
\item [(1)] $\gamma=1$: $h(\cdot)$ reduces to a linear function
and the $p(x)=0$ always has a one positive root, say $p>0$.
$H(\cdot)$ function is convex on $[0, F(p))$ and concave on $(F(p),
+\infty)$.  Hence we have a connected continuation region $(0,
b^*)$.
\item[(2)] $\gamma<1$:  We observe that $\lim_{x\downarrow
0}p(x)=-\infty$, $\lim_{x\uparrow +\infty}p(x)=-\infty$,
$\lim_{x\downarrow 0}p'(x)=+\infty$, and $\lim_{x\uparrow
+\infty}p(x)=0-$.   A direct analysis of $p^{'}(x)$ shows that there
is only one stationary point in $(0, \infty)$ and the number of the
roots of $p(x)=0$ is either $0, 1$ or $2$. Hence in the first two
cases, $H(\cdot)$ is concave on $[0, \infty)$ and the continuation
region is connected.  In the last case where there are two roots,
say $0<p_1<p_2$.  The $H(\cdot)$ function is then concave on $[0,
F(p_1))\bigcup (F(p_2), +\infty)$ and is convex on $(F(p_1),
F(p_2))$. Since $H(\cdot)$ increases and concave on $y\in(F(p_2),
\infty)$, we can conclude that the continuation region is connected
in this case as well.
\end{itemize}
Let us move on to finding an optimal continuation region.  By fixing
$a>0$, let us define $r(x, a)=0$ if $x=0$, $-K$ if $0< x< a$ and
                       $k(x-a)^\gamma-K$ if
                       $x\geq a$.
When we solve (\ref{eq:simplesystem}) for this $a$, it is not easy
(at least analytically) to solve $F^{-1}(y)$ explicitly.  We can
bypass this difficulty by using the first identity of
(\ref{eq:devH}) so that (\ref{eq:simplesystem-absorbing}) with $D=0$
becomes
\begin{equation}\label{eq:bypass}
\left(\frac{r(b, a)}{\varphi(b)}\right)'=\frac{r(b,
a)}{\varphi(b)(F(b)-F(0))-\varphi(a)(F(a)-F(0))}\left(F'(b)+\frac{\varphi'(b)\varphi(a)}{\varphi(b)^2}(F(a)-F(0)\right)
\end{equation}
As in the previous examples, $W_a(\cdot)$ is a straight line passing
$(F(0), 0)$ in the form of $W_a(y)=\beta(y-F(0))$. The value
function $v(x)$ in $x\in(0, b^*)$ is
\begin{align*}
\hat{v}(x)&=\varphi(x)W(F(x))=\beta(F(x)-F(0))\varphi(x)\\
&=\beta^*(\psi(x)-F(0)\varphi(x))=\beta^*e^{\frac{\delta}{2}\frac{(x-m)^2}{\sigma^2}}\left\{\mathcal{D}_{-\alpha/\delta}\left(-\left(\frac{x-m}{\sigma}\right)\sqrt{2\delta}\right)
-F(0)\mathcal{D}_{-\alpha/\delta}\left(\left(\frac{x-m}{\sigma}\right)\sqrt{2\delta}\right)\right\}.
\end{align*}
Therefore, the solution to the problem is
\begin{eqnarray*}
v(x)&=& \begin{cases}
                    \hat{v}(x), &0\leq x \leq b^*, \\
                    \hat{v}(a^*)+k(x-a^*)^\gamma-K, &b^*\leq x.
        \end{cases}
\end{eqnarray*}
This solves the problem. See Figure 2-(b) for the value function in
case of parameters $\delta=0.1$, $m=0.9$, $\sigma=0.35$,
$\alpha=0.105$ for the diffusions.  As for the reward/cost function
parameters, $k=0.7$, $K=0.1$ and $\gamma=0.75$.  The solution is $(a*, b*, \beta)=(0.2192, 0.6220, 0.5749)$.
\end{example}

%----------------------------
% Example Sin
%============================
\begin{example}\normalfont\label{ex:sin}
We show a simple example where we have multiple continuation
regions, the first case of Remark \ref{rmk:pathology}.  Let the
uncontrolled process is a standard Brownian motion $B_t$ and let
$\alpha=0$, $f=0$ and
\begin{equation*}
K(x, y)=-c(\sin x-\sin y)-\delta
\end{equation*}
with $c\in\mathbb{R_+}$ and $\delta\in\mathbb{R_+}$ being some
constant parameters.  We want to solve
\begin{equation*}
v(x)=\sup_{\nu\in S}\ME\left[\sum_{T_i<\tau_0}(\xi_i-\delta)\right].
\end{equation*}
In this case $F(x)=x$ and let us define
\begin{eqnarray*} \nonumber
   R(x, a)=r(x, a) &=&  \begin{cases}
                 0, &x=0, \\
                 -c(\sin x-\sin a)-\delta, &x> 0.
     \end{cases}
\end{eqnarray*}
By solving (\ref{eq:simplesystem}) with some parameter $(c,
\delta)=(10, 0.35)$, we find that $a_k^*=2.75+4k\pi$ and
$b_k^*=3.52+4k\pi$ with $k=0, 1, 2...$.  For all these pairs,
$\beta^*$ has a common value of $9.30$.  Hence all these pairs are
optimal.  This implies that if the initial state $x\in(b_k^*,
b_{k+1}^*)$, then we let the process move until it reaches $b_k^*$
or $b_{k+1}^*$.  If it reaches $b_k^*$ first, then an intervention
is made to $a_k^*$. Now we are in the interval $(b_{k-1}^*,
b_{k}^*)$.  We continue until the process is absorbed at $x=0$.
\end{example}

\section{Conclusions}\label{sec:conclusion}
Before we conclude this article, we shall mention an immediate
extension to two boundary impulse control problems:
\begin{equation}\label{eq:reflect}
    J^\nu(x)=\ME\left[\int_0^\infty e^{-\alpha s}f(X_s)ds+\sum_{i=1}
    e^{-\alpha T_i}C_1(X_{T_i-}, X_{T_i})+\sum_{j=1}
    e^{-\alpha S_j}C_2(X_{S_j-}, X_{S_j})\right]
\end{equation}
and
\begin{equation}\label{eq:reflect2}
    v(x)=\sup_{\nu}J^\nu(x)=J^{\nu^*}(x)
\end{equation}
for all $x\in \mathbb{R}$, where
\begin{equation*}
\nu=(T_1, T_2,....; \zeta_1, \zeta_2,....; S_1, S_2,....;\eta_1,
\eta_2,.....)
\end{equation*}
with $\zeta_i>0$ corresponds to interventions at the upper boundary
at intervention time $T_i$ and $\eta_j<0$ at the lower boundary at
intervention time $S_j$.

Examples of this type include the storage model analyzed by Harrison
et al.~\cite{HST1983} and foreign exchange rate model studied by
Jeanblanc-Picqu\'{e}~\cite{MJ1993}.  The former problem, for
example, is that a controller continuously monitors the inventory so
that the inventory level will not fall below the zero level. He is
allowed to make interventions by increasing and decreasing the
inventory by paying costs associated with interventions.  In this
case, the process remains within some band(s). In other words, the
optimal intervention times are characterized as exit times from an
interval in the form of $(p^*, b^*)$ for $0\leq p^*<b^*$. See Korn
\cite{KO1999} for survey. %(In fact, the other combinations of
%admissible strategies, for example, $\zeta_i>0$ and $\eta_j>0$, will
%lead to a one boundary problem.) This performance measure
%(\ref{eq:reflect}) can be interpreted as follows: $f$ is the
%reward/cost continuously incurred throughout the time. On each
%$T_i\geq 0$, by interventions, one pays cost in the amount of
%$C_1(X_{T_i-}, X_{T_i})$ and the process moves from $X_{T_i-}$ to
%$X_{T_i}$.  On the other hand, on each $S_j \geq 0$, by
%interventions, one pays cost in the amount of $C_2(X_{S_j-},
%X_{S_j})$.

Under suitable assumptions, we can develop a similar argument to the
previous chapters.  Among others, the intervention times can be
characterized as exit times from an interval $(p^*, b^*)$.  We can
also simplify the performance measure,

\begin{multline*}
J^\nu(x)=\ME[1_{\{T_1<S_1\}}e^{-\alpha
T_1}\{C_1(X_{T_1-},X_{T_1})-g(X_{T_1-})+J^\nu(X_{T_1})\}]\\
+\ME[1_{\{T_1>S_1\}}e^{-\alpha
S_1}\{C_2(X_{S_1-},X_{S_1})-g(X_{S_1-})+J^\nu(X_{S_1})\}]+g(x)
\end{multline*}
where $g(x)=\ME\int_0^\infty e^{-\alpha s}f(X_s^0)ds$ as usual.
Again, the last term does not depend on controls, we define $u(x)$
as $u(x)=J^\nu(x)-g(x)$,

\begin{equation}
u(x)=
 \ME[1_{\{\tau_b<\tau_p\}}e^{-\alpha \tau_b}u(b)]
+\ME[1_{\{\tau_b>\tau_p\}}e^{-\alpha \tau_p}u(p)], \quad x\in[p,
b]
%K(x, a)-g(x)+g(a)+u(a), &x\in[b,  \infty)
%\end{cases}
\end{equation}
where $T_1=\tau_b$ and $S_1=\tau_p$ and it follows that
\begin{align}
\frac{u(x)}{\varphi(x)}&=\frac{u(b)(F(x)-F(p))}{\varphi(b)(F(b)-F(p))}+\frac{u(p)(F(b)-F(x))}{\varphi(p)(F(b)-F(p))},
\quad x\in[p, b].\label{eq:complex}
\end{align}
Hence if we define $W\triangleq\frac{u}{\varphi}\circ F^{-1}$, we
have linear characterization again in the transformed space;
\begin{equation}\label{eq:W3}
   W(F(x))=W(F(b))\frac{F(x)-F(p)}{F(b)-F(p)}+W(F(p))\frac{F(b)-F(x)}{F(b)-F(p)}, \quad x\in[p, b].
\end{equation}
and the solution to the problem is described as
\begin{eqnarray*}
u(x)&= \begin{cases}
                        \bar{C}_2(x, q)+u_0(q), & x \leq
                       p\\
                        u_0(x)\triangleq\ME[1_{\{\tau_b<\tau_p\}}e^{-\alpha \tau_b}u(b)]
+\ME[1_{\{\tau_b>\tau_p\}}e^{-\alpha \tau_p}u(p)], & p\leq x \leq b\\
                        \bar{C}_1(x, a)+u_0(a), &b\leq x
        \end{cases}
\end{eqnarray*}
where $\bar{C}_i(x, y)=C_i(x, y)-g(x)+g(y)$ for $i=1$ and $2$.\\

We have studied impulse control problems.  The intervention times
are characterized as exit times of the process from a finite union
of disjoint intervals on the real line. A sufficient condition is
given for the connectedness of the continuation region.  The value
function is shown to be linear in certain transformed space and a
direct calculation method is described for it.  This method can
handle impulse control problems with non-smooth reward and cost
functions.  The finiteness of the value function is shown to be
equivalent to the existence of a concave majorant of a suitable
transformation.  The latter is easier to check by using geometric
arguments.

The new characterization of the value function and optimal
strategies can be extended to other optimization problems, such as
optimal switching, singular stochastic control and combined problems
of optimal stopping and impulse control.  If an optimal strategy
exists in the class of exit times, then the problem can be reduced
to a sequence of optimal stopping problems and an effective
characterization of the value function is possible.

\section{Appendix}
\subsection{Proof of Lemma \ref{lem:Davis}}\label{apx:A} To make
the proof more intuitive, we will work with (\ref{eq:wrec}) rather
than with (\ref{eq:wrec2}) where the integration part is converted
to $g$ functions.  For this purpose, it is convenient to define
the following two operators $\mathcal{M}_o :
\mathcal{H}\rightarrow\mathcal{H}$ and $\mathcal{L}_o
:\mathcal{H}\rightarrow\mathcal{H}$:
\begin{equation}\label{eq:Mo}
\mathcal{M}_ou(x)=\sup_{y\in \mathbb{R}}[K(x, y)+u(y)]
\end{equation}
and
\begin{equation}\label{eq:Lo}
\mathcal{L}_ou(x)=\sup_{\tau\in\mathcal{S}}\ME\left[\int_0^\tau
e^{-\alpha s}f(X_s)ds +e^{-\alpha
\tau}\mathcal{M}_ou(X_{\tau-})\right].
\end{equation}
Hence we can proceed with the arguments developed in Davis
\cite{DV1992}. In terms of the two operators just defined,
(\ref{eq:wrec}) becomes
\begin{align} \label{eq:wn+1}
w_{n+1}(x)&=\sup_{\tau\in\mathcal{S}}\ME\left[\int_0^\tau
e^{-\alpha
s}f(X_s)ds+e^{-\alpha \tau}\mathcal{M}_o w_n(X_{\tau-})\right]\\
&=\mathcal{L}_ow_n(x).
\end{align}
(1) $w_n=v_n$ for all $n$:  Let us now prove $v_n(x)=w_n(x)$ for
all $n\in\mathbb{N}$.  We show already $v_1(x)=w_1(x)$. We assume,
to make an induction argument, that $v_n(x)=w_n(x)$ and prove
$v_{n+1}(x)=w_{n+1}(x)$. We should note that, for each $n$, the
optimization problem in (\ref{eq:Lo}) is an optimal stopping
problem.  Hence by Proposition \ref{prop:4}, we can confine the
set of strategy $S_n$ in (\ref{eq:vn}) into a smaller set, i.e.
barrier strategies;
\begin{equation} \label{eq:barS}
\bar{S}_n\triangleq\{\nu\in S_n: T_i, i\in\mathbb{N}\quad \text{is
an exit time from some interval.} \}.
\end{equation}
Now we proceed with the planned induction argument,
\begin{align*}
v_{n+1}(x)&=\sup_{\nu\in
S_{n+1}}\ME\left[\int_0^{\infty}e^{-\alpha
s}f(X_s)ds+\sum_{T_i}e^{-\alpha T_i}K(X_{T_i-},X_{T_i})\right]\\
&=\sup_{(\tau, \xi) \in \bar{S}_1}\ME\left[\int_0^{\tau}e^{-\alpha
s}f(X_s)ds+e^{-\alpha
\tau}(K(X_{\tau-}, X_{\tau})+v_n(X_\tau))\right]\\
&=\sup_{(\tau, \xi) \in \bar{S}_1}\ME\left[\int_0^{\tau}e^{-\alpha
s}f(X_s)ds+e^{-\alpha
\tau}(K(X_{\tau-}, X_{\tau})+w_n(X_\tau))\right]\\
&=\sup_{\tau \in \mathcal{S}}\ME\left[\int_0^{\tau}e^{-\alpha
s}f(X_s)ds+e^{-\alpha \tau}\mathcal{M}_o
w_n(X_\tau)\right]=\mathcal{L}_ow_n(x)=w_{n+1}(x)
\end{align*}
for all $x\in\mathbb{R}$.  The second equality is due to the
strong Markov property justified by (\ref{eq:barS}). The third
equality is by the induction hypothesis.  This proves the first
statement of the lemma.

(2) $v(x)=\lim_{n\rightarrow \infty}w_n(x)$:  Since $w_n$ is
monotone increasing, the limit $w(x)=\lim_{n\rightarrow}w_n(x)$
exists. Since $S_n\subset S$, $w_n(x)\leq v(x)$. Hence $w(x)\leq
v(x)$.  To show the reverse inequality, we define $S^*$ be a set
of interventions such that
\begin{equation*}
S^*=\{\nu\in\S: J^\nu(x)<\infty \quad \text{for all}\quad
x\in\mathbb{R} \}.
\end{equation*}
Let us assume that $v(x)<+\infty$ and consider strategy $\nu^*\in
S^*$ and another strategy $\nu_n$ that coincides with $\nu^*$ up
to and including time $T_n$ and then takes no further
interventions.
\begin{align*}
J^{\nu^*}(x)-J^{\nu_n}(x)&=\ME\left[\int_{T_n}^\infty e^{-\alpha
s} (f(X_s)-f(X^0_{s}))ds +\sum_{i\geq n+1}e^{-\alpha
T_i}K(X_{T_i-}, X_{T_i}) \right],
\end{align*}
which implies
\begin{equation*}
|J^{\nu^*}(x)-J^{\nu_n}(x)|\leq
\ME\left[\frac{2\|f\|}{\alpha}e^{-\alpha T_n}+\sum_{i\geq
n+1}e^{-\alpha T_i}K(X_{T_i-}, X_{T_i})\right].
\end{equation*}
As $n\rightarrow +\infty$, the right hand side goes to $0$ by the
dominated convergence theorem. Hence it is shown
\begin{equation*}
v(x)=\sup_{\nu\in S^*}J^\nu(x)=\sup_{\nu\in \bigcup_n
S_n}J^\nu(x),
\end{equation*}
so that $v(x)\leq w(x)$.  Next, consider the case of
$v(x)=+\infty$. By Proposition \ref{prop:genK}, we have
$q=+\infty$ in this case. Then by the recursive method described
in Section \ref{sec:first}, we see that $v_1(x)=w_1(x)=\infty$. By
the first statement of this lemma, we can conclude
$v_n(x)=w_n(x)=\infty$ for all $n\in\mathbb{N}$, obtaining
$v(x)=\lim_{n\rightarrow\infty}w_n(x)$. This completes the proof
of the second statement.

(3) $w=\mathcal{L}_ow$ : Since $w_n \uparrow w$, we have the
following chain of equalities:
\begin{align*}
\mathcal{M}_o w(x)&=\sup_{y \in \mathbb{R}}[K(x, y)+w(y)]
=\sup_{y\in \mathbb{R}}\sup_{n\in \mathbb{N}}[K(x, y) +w_n(y)]\\
&=\sup_{n\in \mathbb{N}}\sup_{y\in \mathbb{R}}[K(x,
y)+w_n(y)]=\sup_{n\in \mathbb{N}}\mathcal{M}_o w_n(x).
\end{align*}
In view of this, if we take the limit on the both sides of
(\ref{eq:wn+1}) as $n\rightarrow \infty$, by the monotone
convergence theorem,
\begin{equation*}
w(x)=\sup_{\tau\in\mathcal{S}}\ME\left[\int_0^\tau e^{-\alpha
s}f(X_s)ds+e^{-\alpha \tau}\mathcal{M}_o w(X_\tau)\right].
\end{equation*}
This shows that $w=\mathcal{L}_ow$.  Suppose $w'(x)$ satisfies
$w'=\mathcal{L}_ow'$ and majorizes $g(x)=v_0(x)$. Then
$w'=\mathcal{L}_ow'\geq \mathcal{L}_ov_0=w_1$.  If we assume that
$w'\geq v_n$, then
\begin{equation*}
w'=\mathcal{L}_o w' \geq \mathcal{L}_ov_n=v_{n+1}=w_{n+1}.
\end{equation*}
By the induction argument, we have $w'\geq w_n$ for all $n$, leading
to $w'\geq \lim_{n\rightarrow \infty}w_n=w$. Thus it shows that $w$
is the smallest solution majorizing $g$ of the functional equation,
$w-g=\mathcal{L}(w-g)$. This completes the third statement of the
lemma.
\small{\bibliography{impulse}}
\end{document}